\documentclass[graybox]{svmult}

\usepackage{mathptmx}       % selects Times Roman as basic font
\usepackage{helvet}         % selects Helvetica as sans-serif font
\usepackage{courier}        % selects Courier as typewriter font
\usepackage{type1cm}        % activate if the above 3 fonts are
\usepackage{makeidx}         % allows index generation
\usepackage{graphicx}        % standard LaTeX graphics tool
\usepackage{multicol}        % used for the two-column index
\usepackage[bottom]{footmisc}% places footnotes at page bottom
\usepackage{amsmath,amssymb}
\usepackage{url}

\usepackage{tikz}

\makeindex             % used for the subject index
\newcommand{\card}{\operatorname{card}}

\newcommand{\G}{\mathcal{G}}
\newcommand{\T}{\mathcal{T}}

\newcommand{\CIH}{\|I_H\|_{\mathcal{L}(V)}}

\newcommand{\ddiv}{\operatorname*{div}}

\newcommand{\kernel}{\operatorname*{Ker}}
\newcommand{\support}{\operatorname*{supp}}

\renewcommand{\O}{\mathcal{O}}
\newcommand{\N}{\mathcal{N}}

\newcommand{\R}{\mathbb{R}}
\newcommand{\C}{\mathbb{C}}

\renewcommand{\k}{\kappa}
\newcommand{\Cstab}{C_{\operatorname*{st}}}
\newcommand{\Cor}{\mathcal{C}}

\begin{document}

\title*{Variational Multiscale Stabilization and the Exponential Decay of Fine-scale Correctors}
\author{Daniel Peterseim\thanks{The author gratefully acknowledges support by the Hausdorff Center for Mathematics Bonn and by Deutsche Forschungsgemeinschaft in the Priority Program 1748 "Reliable simulation techniques in solid mechanics. Development of non-standard discretization methods, mechanical and mathematical analysis" under the project
"Adaptive isogeometric modeling of propagating strong discontinuities in heterogeneous materials".}}
\institute{Daniel Peterseim \at Rheinische Friedrich-Wilhelms-Universit{\"a}t Bonn, Institut f{\"u}r Numerische Simulation, Wegelerstr. 6, D-53115 Bonn, Germany, \email{peterseim@ins.uni-bonn.de}}
\maketitle

\abstract{This paper reviews the variational multiscale stabilization of standard finite element methods for linear partial differential equations that exhibit multiscale features. The stabilization is of Petrov-Galerkin type with a standard finite element trial space and a problem-dependent test space based on pre-computed fine-scale correctors. The exponential decay of these correctors and their localisation to local cell problems is rigorously justified. The stabilization eliminates scale-dependent pre-asymptotic effects as they appear for standard finite element discretizations of highly oscillatory problems, e.g., the poor $L^2$ approximation in homogenization problems or the pollution effect in high-frequency acoustic scattering.}

\section{Introduction}
In the past decades, the numerical analysis of partial differential equations (PDEs) was merely focused on the numerical approximation of sufficiently smooth solutions in the asymptotic regime of convergence. In the context of multiscale problems (and beyond), such results have only limited impact because the numerical approximation will hardly ever reach the asymptotic idealised regime under realistic conditions. Although a method performs well for sufficiently fine meshes it may fail completely on coarser (and feasible) scales of discretization. This happens for instance if the PDE exhibits rough and highly oscillatory solutions. Among the prominent applications are the numerical homogenization of elliptic boundary value problems with highly varying non-smooth diffusion coefficient, high-frequency time-harmonic acoustic wave propagation, and singularly perturbed problems such as convection-dominated flow. 

The numerical approximation of such problems by finite element methods (FEMs) or related schemes is by no means straight-forward. The pure approximation (e.g. interpolation) of the unknown solutions by finite elements already requires high spatial resolution to capture fast oscillations and heterogeneities on microscopic scales. 
When the function is described only implicitly as the solution of some partial differential equation, its approximation faces further scale-dependent pre-asymptotic effects caused by the under-resolution of relevant microscopic data. Examples are the poor $L^2$ approximation in homogenization problems (see Fig.~\ref{fig:numhom}) and the pollution effect \cite{BabSau} for Helmholtz problems with large wave numbers (see Fig.~\ref{fig:helmholtz}). We shall emphasise that, in the latter case, the existence and uniqueness of numerical approximations may not even be guaranteed in pre-asymptotic regimes.
\definecolor{matlabred}{RGB}{216,82,24}
\begin{figure}[b]
%\sidecaption[b]
% For example, with the graphicx style use
\includegraphics[width=\textwidth]{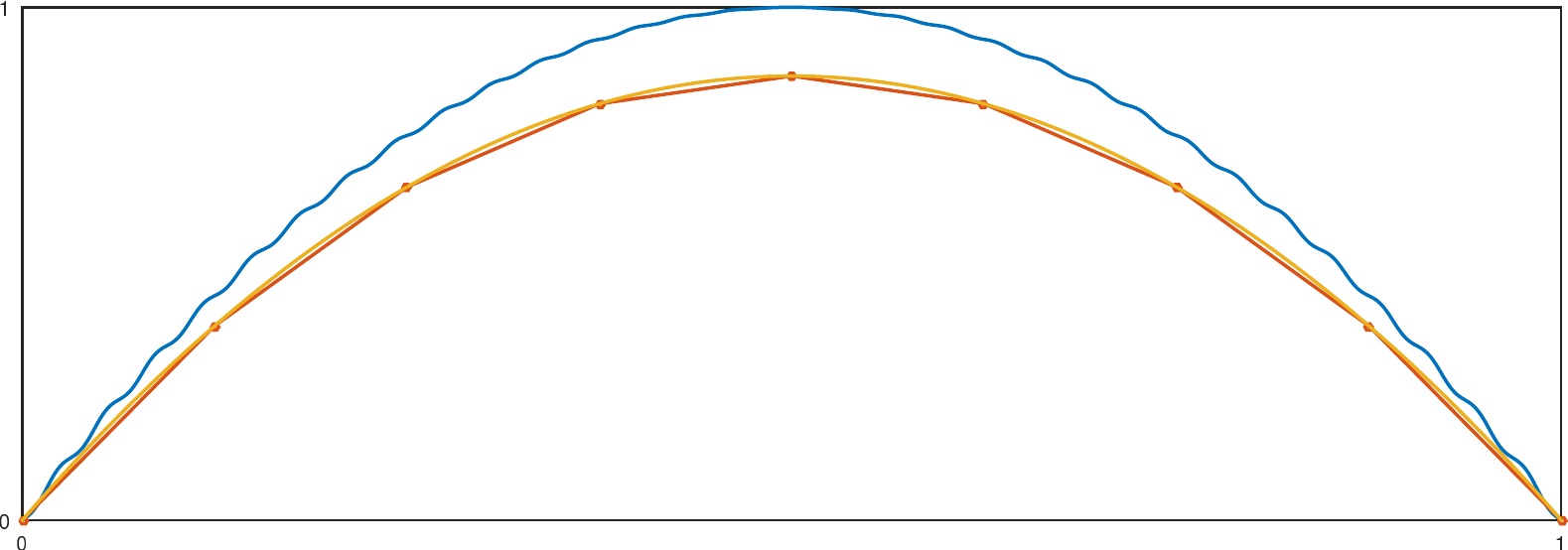}
%
% If no graphics program available, insert a blank space i.e. use
%\picplace{5cm}{2cm} % Give the correct figure height and width in cm
%
\caption{Failure of FEM in homogenization problems: Consider the periodic problem $-\tfrac{d}{dx}A_\varepsilon(x)\tfrac{d}{dx}u_\varepsilon(x) = 1$ in the unit interval with homogeneous Dirichlet boundary condition, where $A_\varepsilon(x):=(2+\cos(2\pi x/\varepsilon))^{-1}$ for some small parameter $\varepsilon>0$. The solution $u_{\varepsilon}=4(x-x^2)-4\varepsilon\left(\frac{1}{4\pi}\sin(2\pi\tfrac{x}{\varepsilon})-\frac{1}{2\pi}x\sin(2\pi\tfrac{x}{\varepsilon})-\frac{\varepsilon}{4\pi^2}\cos(2\pi\tfrac{x}{\varepsilon})+\frac{\varepsilon}{4\pi^2}\right)$ is depicted in blue for $\varepsilon=2^{-5}$. The P1-FE approximation (\textcolor{matlabred}{$\circ$}) on a uniform mesh of width $h$ interpolates the curve $x\mapsto 2\sqrt{3}(x-x^2)$ whenever $h$ is some multiple of the characteristic length scale $\varepsilon$ and, hence, fails to approximate $u_\varepsilon$ in any reasonable norm in the regime $h\geq\varepsilon$.}
\label{fig:numhom}       % Give a unique label
\end{figure}

\begin{figure}[b]
%\sidecaption[b]
% For example, with the graphicx style use
\includegraphics[width=\textwidth]{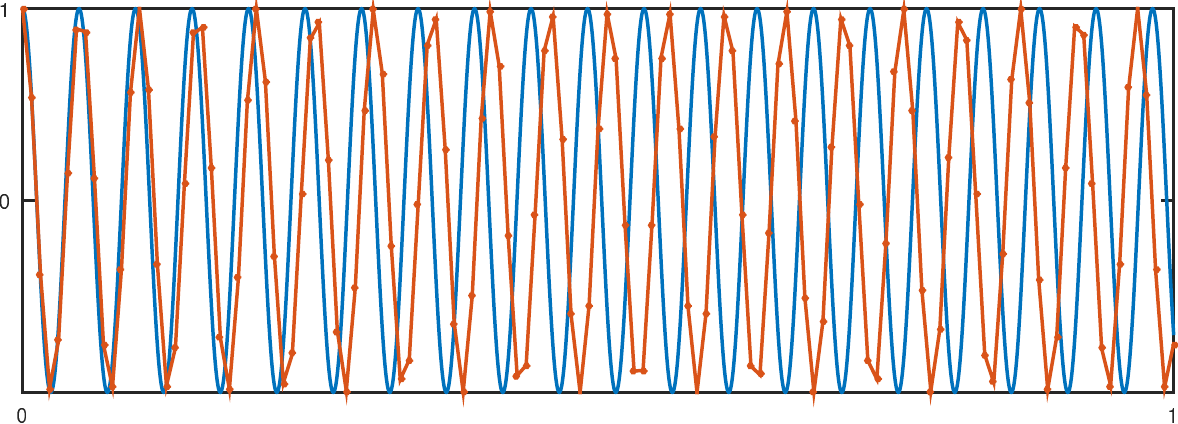}
%
% If no graphics program available, insert a blank space i.e. use
%\picplace{5cm}{2cm} % Give the correct figure height and width in cm
%
\caption{Numerical dispersion in Helmholtz problems: Consider $-\tfrac{d^2}{dx^2}u_\varepsilon(x) - \kappa^2 u(x)= 0$ in the unit interval with $u(0)=1$ and $\tfrac{d}{dx}u(1)=-i\kappa u(1)$ for some large parameter $\kappa>0$. The solution $u_{\kappa}=\exp(-i\kappa x)$ is depicted in blue for $\kappa=2^{7}$. The P1-FE approximation (\textcolor{matlabred}{$\boldsymbol{\circ}$}) on a uniform mesh of width $h=2^{-7}>6\cdot(\text{wave length})$ fails to approximate $u_\kappa$ due to the accumulation of phase errors.}
\label{fig:helmholtz}       % Give a unique label
\end{figure}

Such situations require the stabilization of standard methods so that eventually a meaningful approximation on reasonably coarse scales of discretization becomes feasible. This paper aims to present a general framework for the stabilization of FEMs for multiscale problems with the aim to significantly reduce or even eliminate pre-asymptotic effects due to under-resolution. Our starting point will be the Variational Multiscale Method (VMS) originally introduced in \cite{Hughes:1995,MR1660141}. The method provides an abstract framework how to incorporate missing fine-scale effects into numerical problems governing coarse-scale behaviour \cite{MR2300286}. One may interpret the VMS as a Petrov-Galerkin method using standard FE trial spaces and an operator-dependent test space that needs to be precomputed in general. 

The construction of this operator-dependent test space is based on some stable projection onto the standard FE trial space and a corresponding scale decomposition of a function into its FE part given by the projection (the macroscopic/coarse-scale part) and a remainder that lies in the kernel of the projection operator (the microscopic/fine-scale part). The test functions are computed via a problem-de\-pen\-dent projection of the trial space into the space of fine-scale functions. This requires the solution of variational problems in the kernel of the projection -- the fine-scale corrector problems. It has been observed empirically in certain applications that the Green's function associated with these fine-scale corrector problems -- the so-called fine-scale Green's function \cite{MR1660141} -- may exhibit favourable exponential decay properties \cite{axel_phd,MR1660141} even though the decay of the classical full scale Green's function is only algebraic. It is this exponential decay property that allows one to turn the VMS into a feasible numerical method \cite{axel_phd,LM07}. 

Very recently, the exponential decay was rigorously proved for the first time in \cite{MP14} in the context of multi-dimensional numerical homogenization. A key ingredient of the proof of \cite{MP14} is the use of a (local) quasi-interpolation operator for the scale decomposition. Although the method of \cite{MP14} still fits into the general framework of the VMS, it uses a different point of view on the method based on the orthogonalisation of coarse and fine scales with respect to the inner product associated with a symmetric and coercive model problem. This is why the method is now referred to as the Localized Orthogonal Decomposition (LOD) method. Subsequent work showed that the ideas of \cite{MP14} can be generalised to other discretization techniques such as discontinuous Galerkin \cite{MR3084163,EGMP13,E15}, Petrov-Galerkin formulations \cite{EGH15}, mixed methods \cite{HHM15} and mesh-free methods \cite{HMP14}. Moreover, the method can also be reinterpreted in terms of the multiscale finite element method with special oversampling \cite{HP13}. The class of problems that have been analysed by now includes semi-linear problems \cite{HMP13}, high-contrast problems \cite{Peterseim.Scheichl:2014,Brown.Peterseim:2014}, rough boundary conditions \cite{HM14}, problems on complicated geometries \cite{ELM15}, linear and non-linear eigenvalue problems \cite{MP12,Henning.Mlqvist.Peterseim:2013,MP15}, parabolic problems \cite{2015arXiv150408140M}, wave propagation \cite{AH14,Peterseim2014,Gallistl.Peterseim:2015} and parametric problems \cite{AH14a}.

This survey aims to reinterpret all those results in the abstract stabilization framework of the original VMS (Section~\ref{s:abstract}) and aims to illustrate how the exponential decay of the fine-scale Green's function can be quantified (Section~\ref{s:localization}). We will show how these abstract results lead to super-localised numerical homogenization \cite{MP14,HP13} (Section~\ref{s:hom}) and pollution-free time-harmonic acoustic scattering (Section~\ref{s:helmholtz}) \cite{Peterseim2014,Gallistl.Peterseim:2015}. Section~\ref{s:final} contains some final remarks and also also identifies methodological similarities and differences with some other numerical approaches that receive great attention these days, e.g., discontinuous Petrov-Galerkin methods (dPG) \cite{MR2743600}, Trefftz-type methods \cite{Gittelson2009} and Isogeometric Analysis (IGA) \cite{BCH09}. 

\section{Abstract variational multiscale stabilization}\label{s:abstract}
This section is concerned with an abstract variational problem in a complex Hilbert space $V$ as it appears for the weak formulation of second order PDEs. In this context, $V$ is typically some closed subspace of the Sobolev space $W^{1,2}(\Omega;\mathbb{C}^m)$ for some bounded Lipschitz domain $\Omega\subset\mathbb{R}^d$. Let $a$ denote a bounded sesquilinear form on $V\times V$ and let $F\in V'$ denote a bounded linear functional on $V$. We wish to find $u\in V$ satisfying the linear variational problem
\begin{equation}\label{e:modelabstract}
\forall v\in V:\; a(u,v) = \overline{F(v)}.
\end{equation}
We assume that the sesquilinear form $a$ satisfies the inf-sup condition
\begin{equation}\label{e:BNB1}
\alpha:=\inf_{0\neq v\in V}\sup_{0\neq w\in V}\frac{a(v,w)}{\|v\|_V\|w\|_V}=\inf_{0\neq w\in V}\sup_{0\neq v\in V}\frac{a(v,w)}{\|v\|_V\|w\|_V}>0.
\end{equation}
Under this condition, the abstract problem \eqref{e:modelabstract} is well-posed, i.e., for all $F\in V'$ there exists a unique solution $u\in V$ and the a priori bound
\begin{equation*}
\|u\|_V\leq \alpha^{-1}\|F\|_{V'}
\end{equation*}
holds true; see, e.g., \cite{B71}. 

We wish to approximate the unknown solution $u$ of \eqref{e:modelabstract} by some computable function. The standard procedure for approximation is the Galerkin method which simply chooses a finite-dimensional subspace $V_H\subset V$ (that contains simple functions such as piecewise polynomials) and restricts the variational problem \eqref{e:modelabstract} to this subspace. Usually, $V_H$ belongs to some family of spaces parametrised by some abstract discretization parameter $H$, for instance the mesh size. This parameter (or set of parameters) provides some control on the approximation properties of $V_H$ as $H\rightarrow 0$ at the price of an increasing computational cost in the sense of $\dim V_H\rightarrow \infty.$
The Galerkin method seeks a function $G_Hu\in V_H$ satisfying
\begin{equation}\label{e:Galerkin}
\forall v_H\in V_H:\; a(G_Hu,v_H) = \overline{F(v_H)}\;\;\bigl(\;=a(u,v_H)\;\bigr).
\end{equation}
Recall that the well-posedness of the original problem \eqref{e:modelabstract} does not imply the well-posedness of the discrete variational problem \eqref{e:Galerkin} but needs to be checked for the particular application via discrete versions of the inf-sup condition \eqref{e:BNB1}. In many cases, such conditions are only satisfied for $H$ sufficiently small. This means that there is some threshold complexity for computing any Galerkin approximation and this threshold can be out of reach. Even if a Galerkin solution $G_Hu$ exists and is computable, it might not provide the desired accuracy or does not reflect the relevant characteristic features of the solution, as we have seen in the introduction.

Therefore, we are interested in computing projections onto the discrete space $V_H$ other than the Galerkin projection. Let $I_H:V\rightarrow V_H$ denote such a linear surjective projection operator and let us assume that it is bounded in the sense of the space $\mathcal{L}(V)$ of linear operators from $V$ to $V$ with finite operator norm 
\begin{equation*}
 \|I_H\|_{\mathcal{L}(V)}:=\sup_{0\neq v\in V}\frac{\|I_Hv\|_V}{\|v\|_V}<\infty.
\end{equation*}
Implicitly, we also assume that this operator norm does not depend on the discretization parameter $H$ in a critical way. Possible choices of $I_H$ include the orthogonal projection onto $V_H$ with respect to the inner product of $V$ or any Hilbert spaces $L\supset V$ containing $V$ and mainly (local) quasi-interpolation operators of Cl\'ement or Scott-Zhang type as they are well-established in the finite element community in the context of a posteriori error estimation \cite{MR0400739,MR1011446,MR1736895,CV}.

\subsection{Petrov-Galerkin characterisation of finite element projections}\label{s:PG}
The Galerkin projection $G_H$ is designed in such a way that its computation requires only the known data $F$ associated with the unknown solution $u$. This section characterises the projection $I_H\in\mathcal{L}(V)$ as a Petrov-Galerkin discretization of \eqref{e:modelabstract} using $V_H$ as the trial space and a non-standard test space $W_H\subset V$ that depends on the problem and the projection. The definition of $W_H$ rests on the trivial observation that, for any $v\in V$, 
\begin{equation}\label{e:test1}
a(I_H u,v) = \overline{F(v)} - a(u-I_H u,v).
\end{equation}
The choice of a test function $v$ in the subspace
\begin{equation}\label{e:test}
W_H:=\left\{w\in V\;\vert\;\forall z\in\kernel I_H:\,a(z,w)=0\right\}
\end{equation}
annihilates the second term on the right-hand side of \eqref{e:test1} and, hence, 
\begin{equation*}
a(I_H u,w_H) = \overline{F(w_H)}
\end{equation*}
holds for all $w_H\in W_H$. This shows that $I_H u$ is a solution of the Petrov-Galerkin method: Find $u_H\in V_H$ such that
\begin{equation}\label{e:PGalerkin}
\forall w_H\in W_H:\; a(u_H,w_H) = \overline{F(w_H)}.
\end{equation}
This characterisation of $I_H$ is well known from the variational multiscale method as it is presented in \cite{MR2300286}. 

The question whether or not \eqref{e:PGalerkin} has a unique solution can not be answered under the general assumptions made so far. We need to assume the missing uniqueness to be able to proceed and one way of doing this is to assume that the dimensions of trial and test space are equal, 
\begin{equation}\label{e:BNBH}
\operatorname{dim}W_H = \operatorname{dim}V_H.
\end{equation}
In the present setting with a bounded operator $I_H$, this condition is equivalent to the well-posedness of the discrete variational problem \eqref{e:PGalerkin}, i.e., it admits a unique solution $u_H=I_H u\in V_H$ and
\begin{equation*}%\label{e:apriori}
\|u_H\|_V\leq \|I_H\|_{\mathcal{L}(V)}\|u\|_V\leq \frac{\CIH}{\alpha}\|F\|_{V'}.
\end{equation*}
The a priori estimate in turn implies a lower bound of the discrete inf-sup constant of the Petrov-Galerkin method by the quotient of the continuous inf-sup constant $\alpha$ and the continuity constant of $I_H$, 
\begin{equation*}\label{e:BNB1H}
\inf_{0\neq v_H\in V_H}\sup_{0\neq w_H\in W_H}\frac{a(v_H,w_H)}{\|v_H\|_V\|w_H\|_V}\geq\frac{\alpha}{\CIH}\leq\inf_{0\neq w_H\in W_H}\sup_{0\neq v_H\in V_H}\frac{a(v_H,w_H)}{\|v_H\|_V\|w_H\|_V}.
\end{equation*}
\smallskip

The test space $W_H$ is the ideal test space for our purposes in the following sense. Assuming that we have access to it, the method \eqref{e:PGalerkin} would enable us to compute $I_H u$ without the explicit knowledge of $u$. Although this will rarely be the case, we will see later that $W_H$ can be approximated very efficiently in relevant cases. The discrete inf-sup conditions then indicate that the sufficiently accurate approximation of $W_H$ will not harm the method, its stability properties or its subsequent error minimisation properties.

The continuity of the projection operator $I_H$ readily implies the quasi-optimality of the Petrov-Galerkin method \eqref{e:PGalerkin},
\begin{equation}\label{e:quasibestV}
 \|u-u_H\|_V=\|(1-I_H)u\|_V\leq \|I_H\|_{\mathcal{L}(V)} \min_{v_H\in V_H}\|u-v_H\|_V.
\end{equation}
Here, we have used that $\|I_H\|_{\mathcal{L}(V)}=\|1-I_H\|_{\mathcal{L}(V)}$; see e.g. \cite{MR2279449}. 
More importantly, the same arguments show that the Petrov-Galerkin method is quasi-optimal with respect to any other Hilbert space $L\supset V$ with norm $\|\cdot\|_L$ whenever $I_H\in\mathcal{L}(L)$,
\begin{equation*}
 \|u-u_H\|_L\leq \|I_H\|_{\mathcal{L}(L)} \min_{v_H\in V_H}\|u-v_H\|_L.
\end{equation*}

These quasi-optimality make the ansatz very appealing and motivates further investigation. Hence, in the remaining part of the paper, it is our aim to turn the method into a feasible numerical scheme while preserving these properties to a large extent. Although the discrete stability of the method depends on the stability properties of the original problem and, hence, on parameters such as the frequency in scattering problems, the quasi-optimality depends only on $I_H$ and not necessarily on the problem.   

\subsection{Characterisation of the ideal test space}\label{s:test}
A practical realisation of the Petrov-Galerkin method \eqref{e:PGalerkin} requires a choice of bases in the discrete trial $V_H$ and test space $W_H$. These choices will have big impact on the computational complexity. The underlying principle of finite elements is the locality of the bases which yields sparse linear systems and offers the possibility of linear computational complexity with respect to the number of degrees of freedom. Let $\{\lambda_j\;\vert\;j=1,2,\ldots,N_H=\operatorname{dim}V_H\}$ be such a local basis of $V_H$.   

We shall derive a basis of the test space $W_H$ defined in \eqref{e:test} by mapping the trial basis onto a test basis via some bijective operator $\mathcal{T}$, a so-called trial-to-test operator. Due to Assumption \eqref{e:BNBH} such an operator exists, but there are many choices and we have to make a design decision. Our choice is  that 
\begin{equation}\label{e:design}
I_H\circ \T = id
\end{equation}
which is consistent with almost all existing practical realisations of the method but one might as well consider 
distance minimisation
$\|(1-\T)v_H\|_V=\min_{w_H\in W_H}\|v_H-w_H\|_V$.

The condition \eqref{e:design} fixes the (macroscopic) finite element part $I_H\T v_H=v_H$ of $\T v_H$ while the fine scale remainder $(1-I_H)\T v_H$ is determined by the variational condition in the definition of $W_H$. Given $v_H\in V_H$, $(1-I_H)\T v_H \in\kernel I_H$ satisfies
\begin{equation}\label{e:correction}
 \forall z\in\kernel I_H:a(z,(1-I_H)\T v_H)=-a(z,v_H).
\end{equation}
This problem is referred to as the fine scale corrector problem for $v_H\in V_H$. Note that $v_H$ can be replaced with any $v\in V$ so that $(1-I_H)\T$ can be understood as an operator from $V$ into $\kernel I_H$. We usually denote this operator the fine scale correction operator and write $\Cor:=(1-I_H)\T$. This operator is the Galerkin projection from $V_H$ (or $V$) into $\kernel I_H$ related to the adjoint of the sesquilinear form $a$. It depends on the underlying variational problem and equips test functions with problem related features that are not present in $V_H$. In the context of elliptic PDEs, $\Cor$ is called the finescale Green's operator \cite{Hughes:1995,MR2300286}. 

For this construction to work we need to assume the well-posedness of the corrector problem \eqref{e:correction}, i.e., there is some constant $\beta>0$ such that
\begin{equation}\label{e:BNB1fine}
\inf_{0\neq v\in \kernel I_H}\sup_{0\neq w\in \kernel I_H}\frac{a(v,w)}{\|v\|_V\|w\|_V}\geq \beta\leq\inf_{0\neq w\in \kernel I_H}\sup_{0\neq  v\in \kernel I_H}\frac{a(v,w)}{\|v\|_V\|w\|_V}.
\end{equation}
As for the Galerkin projection $G_H$ onto $V_H$, these inf-sup conditions do not follow from their continuous counterparts \eqref{e:BNB1} (unless $a$ is coercive) and they might hold for sufficiently small $H$ only. However, we were able to show in the context of the Helmholtz model problem of Section~\ref{s:helmholtz} that \eqref{e:BNB1fine} holds in a much larger regime of the discretization parameter $H$ than the corresponding conditions for the standard FEM do. In any case,  condition~\eqref{e:BNB1fine} implies that the trial-to-test operator $\T=1+\Cor$ is a bounded linear projection operator from $V$ to $W_H$ with operator norm
\begin{equation*}
 \|\T\|_{\mathcal{L}(V)}=\|1-\T\|_{\mathcal{L}(V)}=\|\Cor\|_{\mathcal{L}(V)}\leq\frac{C_a}{\beta},
\end{equation*}
where $C_a$ denotes the continuity constant of the sesquilinear form $a$. 
Moreover, $\T\vert_{V_H}:V_H\rightarrow W_H$ is invertible with $(\T\vert_{V_H})^{-1}=I_H$ and $\{\T\lambda_j\;\vert\;j=1,2,\ldots,N_H\}$ defines a basis of $W_H$ with 
\begin{equation*}
 \frac{1}{\CIH}\|\lambda_j\|_V\leq \|\T\lambda_j\|_V \leq \tfrac{C_a}{\beta}\|\lambda_j\|_V,\quad 1\leq j \leq N_H.
\end{equation*}

In general, it cannot be expected (apart from one-dimensional exceptions where $\kernel I_H$ is a broken Sobolev space \cite{MR2300286}) that the $\T\lambda_j$ have local support. On the contrary, their support will usually be global. However, we will show in the next section that they decay very fast in relevant applications; for illustrations see Section~\ref{s:hom}.

An important special case of the model problem \eqref{e:modelabstract} is the hermitian case. Note that hermiticity is preserved by the Petrov-Galerkin method in the following sense. For any $u_H,v_H\in V_H$, it holds that 
\begin{equation*}
a(u_H,\T v_H)=a(\T u_H,\T v_H) =\overline{a(\T v_H,\T u_H)}=\overline{a(v_H,\T u_H)}.
\end{equation*}
However, this hermiticity is typically lost once $\T$ is replaced with some approximation $\T_\ell$. In order to avoid a lack of hermiticity, previous papers such as \cite{MP14} mainly used a variant of the method with $W_H$ as the test and trial space. If hermiticity is important, one should follow this line. In this paper, we trade hermiticity for a cheaper method that avoids any costly communication between the fine-scale correctors that would be necessary in the hermitian version. 

If the problem is non-hermitian, one might still consider a modified trial space based on the adjoint of $\T$ to improve approximation properties; see \cite{LM09,MRX,Peterseim2014} for details.     
In a setting with a modified trial space, further generalisations are possible. Since $V_H$ does not appear any more in the method, its conformity can be relaxed as it was recently proposed in \cite{2015arXiv150303467O} in the context of a multilevel solver for Poisson-type problems with $L^\infty$  coefficients. This approach enables one to compute very general quantities of the solution such as piecewise mean values.

\section{Exponential decay of fine-scale correctors}\label{s:localization}
In many cases, the fine-scale correctors (i.e. the solutions of the fine-scale corrector problems \eqref{e:correction}) have decay properties better than those of the Green's function associated with the underlying full-scale partial differential operator. To elaborate on this, we shall now assume that the space $V$ is a closed subspace of $W^{1,2}(\Omega)$ with a local norm (the notation $\|\cdot\|_{V,\omega}$ means that the $V$-norm is restricted to some subdomain $\omega\subset\Omega$). Moreover, the sesquilinear form $a$ is assumed to be local. This is the natural setting for scalar second order PDEs. The subsequent arguments can be easily generalised to vector-valued problems. 

To be more precise regarding the locality of the basis mentioned above, we shall associate the basis functions of $V_H$ with a set of geometric entities $\N_H$ called nodes (e.g. the vertices of a triangulation) and assume that these nodes are well distributed in the domain $\Omega$ in the sense of local quasi-uniformity. In this context, $H$ refers to the maximal distance between nearest neighbours (the mesh size).  
Given some node $z\in\N_H$ and the corresponding basis function $\lambda_z\in V_H$, set the corrector $\phi_z=\Cor\lambda_z$ and recall from \eqref{e:correction} that 
\begin{equation*}
 a(w,\phi_z)=-a(w,\lambda_z)
\end{equation*}
for all $w\in \kernel I_H$. 

We aim to show that there are constants $c>0$ and $C>0$ independent of $H$ and $R$ such that 
\begin{equation}\label{e:decay0}
 \|\Cor \lambda_z\|_{V,\Omega\setminus B_{R}(z)}=\|\phi_z\|_{V,\Omega\setminus B_{R}(z)}\leq C\exp\left(-c\frac{R}{H}\right)\|\Cor \lambda_z\|_V,
\end{equation}
where $B_R(z)$ denotes the ball of radius $R>0$ centred at $z$. 

We shall show how this result can be established and what kind of assumptions have to be made. Let $R>2H$ and $r:=R-H>H$ and let $\eta\in W^{1,\infty}(\Omega;[0,1])$ be some cut-off function with $\eta=0$ in $\Omega\setminus B_R(z)$, $\eta=1$ in $B_r(z)$, and 
\begin{equation}\label{e:etaestimate}
\|\nabla\eta\|_{L^\infty(\Omega)}\leq C_\eta H^{-1}
\end{equation}
for some generic constant $C_\eta$. 
In general, the fine-scale space $\kernel I_H$ is not closed under multiplication by a cut-off function and we will need to project the truncated function $\eta\phi_z$ back into $\kernel I_H$ by the operator $1-I_H$. We assume that the concatenation of multiplication by $\eta$ and $(1-I_H)$ is stable and quasi-local in the sense that 
\begin{equation}\label{e:cutclos}
 \forall w\in \kernel I_H:\;\|(1-I_H)(\eta w)\|_{V,B_R(z)\setminus B_r(z)}\leq C_{\eta,I_H} \|w\|_{V,B_{R'}(z)\setminus B_{r'(z)}}
\end{equation}
holds with $r':=r-mH $ and $R':=R+mH$ and generic constants $C_{\eta,I_H}>0$ and $m\in\mathbb{N}_0$ independent of $H$ and $z$. Although 
the multiplication by $\eta$ is not a stable operation in the full space $V$ (think of a constant function), this result is possible in the space of fine scales for example if $I_H$ enjoys quasi-local stability and approximation properties; see Section~\ref{s:hom} below for an example. The quasi-locality of $I_H$ is also used in the next argument. 

Assuming that the inf-sup condition \eqref{e:BNB1fine} holds, the corrector $\phi_z$ satisfies 
\begin{eqnarray*}
 \|\phi_z\|_{V,\Omega\setminus B_R(z)}&=&\|(1-I_H)\phi_z\|_{V,\Omega\setminus B_R(z)}\\
 &\leq& \|(1-I_H)((1-\eta)\phi_z)\|_{V}\\
 &\leq& \beta^{-1} a(w,(1-I_H)((1-\eta)\phi_z))\\
 &=& \beta^{-1} \left(a(w,\phi_z)-a(w,(1-I_H)(\eta\phi_z))\right)
\end{eqnarray*}
for some $w\in \kernel I_H$ with $\|w\|_V=1$. Since $\support \left((1-I_H)((1-\eta)\phi_z)\right)\subset\Omega\setminus B_r(z)$ there is a good chance to actually find a function $w$ with 
\begin{equation*}
\support w\subset \support \left((1-I_H)((1-\eta)\phi_z)\right)\subset\Omega\setminus B_r(z).  
\end{equation*}
Of course, this is an assumption that needs to be verified in the particular application. Under this condition, the term $a(w,\phi_z)=a(w,\lambda_z)$ vanishes because the supports of $w$ and $\lambda_z$ have no overlap. This and \eqref{e:cutclos} imply
\begin{eqnarray*}
 \|\phi_z\|_{V,\Omega\setminus B_R(z)}&\leq& \beta^{-1}C_aC_{\eta,I_H}\|\phi_z\|_{V,B_{R'}(z)\setminus B_{r'}(z)}\\
 &=& \beta^{-1}C_aC_{\eta,I_H}\left(\|\phi_z\|_{V,\Omega\setminus B_{r'}(z)}^2-\|\phi_z\|_{V,\Omega\setminus B_{R'}(z)}^2\right)^{1/2},
\end{eqnarray*}
where $C_a$ denotes the continuity constant of the sesquilinear form $a$. 
Hence, the contraction
\begin{equation*}
 \|\phi_z\|_{V,\Omega\setminus B_{R'}(z)}^2\leq \frac{C'}{1+C'}\|\phi_z\|_{V,\Omega\setminus B_{R'-(2m+1)H}(z)}^2
\end{equation*}
holds with $C':=(\beta^{-1}C_aC_{\eta,I_H})^2$. The iterative application of this estimate with $R'\mapsto R'-(2m+1)H$ plus relabelling $R'\mapsto R$ leads to 
the conjectured decay result \eqref{e:decay0} with constants
$C:=(\frac{C'}{1+C'})^{-\frac{1}{2(2m+1)}}$ and $c:=\left|\log(\frac{C'}{1+C'})\right|\frac{(1)}{2(2m+1)}>0$.
\medskip

The exponential decay motivates and justifies the localisation of the fine-scale corrector problems to local subdomains of diameter $\ell H$ where $\ell\in\mathbb{N}$ is a new discretization parameter, the so-called oversampling parameter. It controls the perturbation with respect to the ideal global correctors. We will explain this localisation procedure on the basis of an example in Section~\ref{s:hom} below. 
As a rule of thumb, the localisation to subdomains of diameter $\ell H$ will introduce an error of order $\O(\exp(-\ell) )$. As long as this error is small when compared with the inf-sup constant $\alpha\|I_H\|_{\mathcal{L}(V)}^{-1}$ of the ideal method, the stability and approximation properties of the method will be largely preserved.  

%This section gave a rather abstract overview on Petrov-Galerkin stabilization in the sense of the variational multiscale method. With the help of two model problems we will show in the subsequent sections that the methodology is indeed useful. 

\section{Application to numerical homogenization of elliptic PDEs}\label{s:hom}
The first prototypical model problem concerns the diffusion problem $
 -\operatorname{div}A\nabla u = f$  
in some bounded domain $\Omega\subset\R^d$ with homogeneous Dirichlet boundary condition. The difficulty is the strongly heterogeneous and highly varying (non-periodic) diffusion coefficient $A$. The heterogeneities and oscillations of the coefficient may appear on several non-separable scales. We assume that the diffusion matrix $A\in L^\infty\left(\Omega,\mathbb{R}_{\mathrm{sym}}^{d\times d}\right)$ is symmetric and uniformly elliptic with
\begin{equation*}%
0<\alpha =\underset{x\in\Omega}{\operatorname{ess}\inf}%
\inf\limits_{v\in\mathbb{R}^{d}\setminus\{0\}}\dfrac{(A( x) v)\cdot v
}{v\cdot v}.
\end{equation*}
Given $f\in L^{2}( \Omega) $, we wish to find the unique weak solution $u\in V:=H_{0}^{1}(\Omega) $
such that%
\begin{equation}\label{e:modelhom}
a\left(  u,v\right):=\int_{\Omega}(A\nabla u)\cdot \nabla
v =\int_{\Omega}fv=:F( v) \quad\text{for all } v\in V.
\end{equation}
It is well known that classical polynomial based FEMs can perform arbitrarily badly for such problems, see e.g.~\cite{MR1648351}. This is due to the fact that finite elements tend to average unresolved scales of the coefficient and the theory of homogenization shows that this way of averaging does not lead to meaningful macroscopic approximations. This is illustrated in the introduction. In the simple periodic example of Fig.~\ref{fig:numhom}, the averaging of the inverse of the diffusion coefficient $A$ (harmonic averaging) would have lead to the correct macroscopic representation. 

In computational homogenization, the impact of unresolved microstructures encoded in the rough coefficient $A$ on the overall process is taken into account by the solution of local microscopic cell problems. While many approaches are empirically successful and robust for certain multiscale problems, the question whether such methods are stable and accurate beyond the strong assumptions of analytical homogenization regarding scale separation or even periodicity remained open for a long time. Only recently, the existence of an optimal approximation of the low-regularity solution space by some arbitrarily coarse generalised finite element space (that represents the homogenised problem) was shown in \cite{BL11} and \cite{ALBasis1}. However, the constructions therein include prohibitively expensive global solutions of the full fine scale problem or the solution of more involved eigenvalue problems. The first efficient and feasible construction, solely based on the solution of localised microscopic cell problems, was given and rigorously justified in \cite{MP14} and later optimised and generalised in \cite{HP13,HMP14}. A different approach with presumably similar properties was later suggested by \cite{BOZ13} along with the notion of sparse super-localisation that reflects the locality of the discrete homogenised operator (similar to the sparsity of standard finite element matrices).

\begin{figure}[tb]
%\sidecaption[t]
% For example, with the graphicx style use
\includegraphics[width=.5\textwidth]{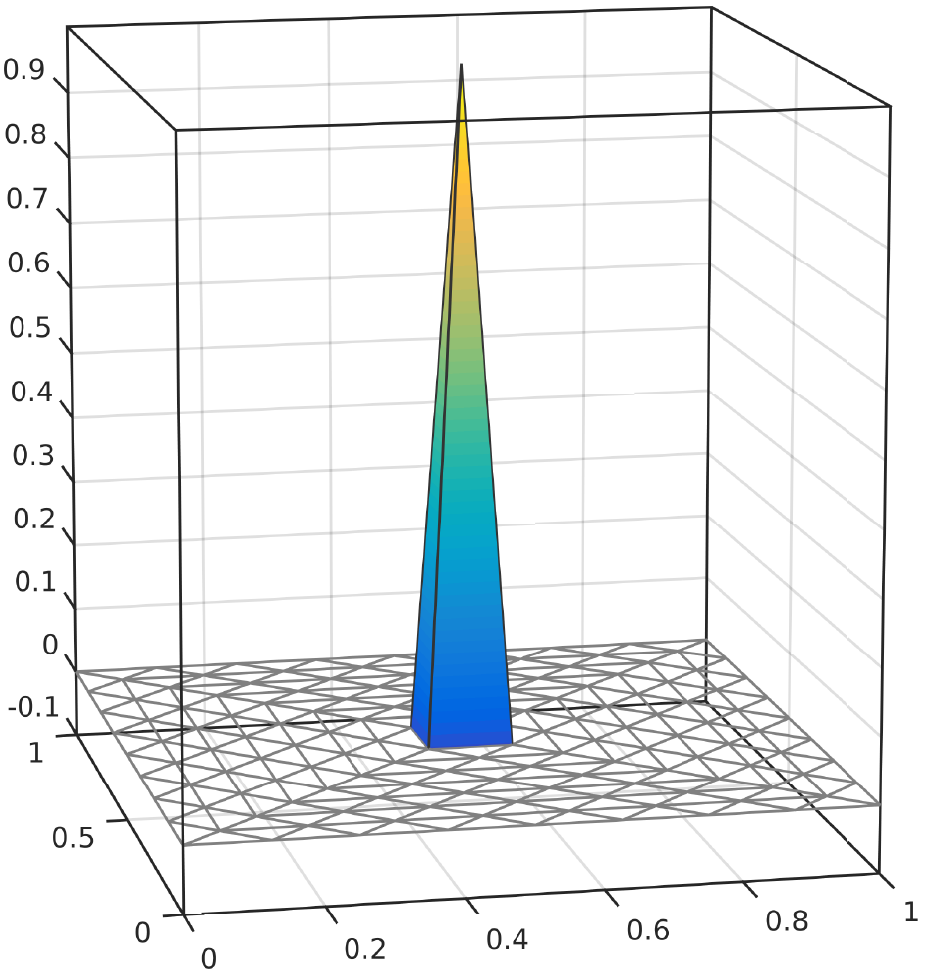}
\includegraphics[width=.5\textwidth]{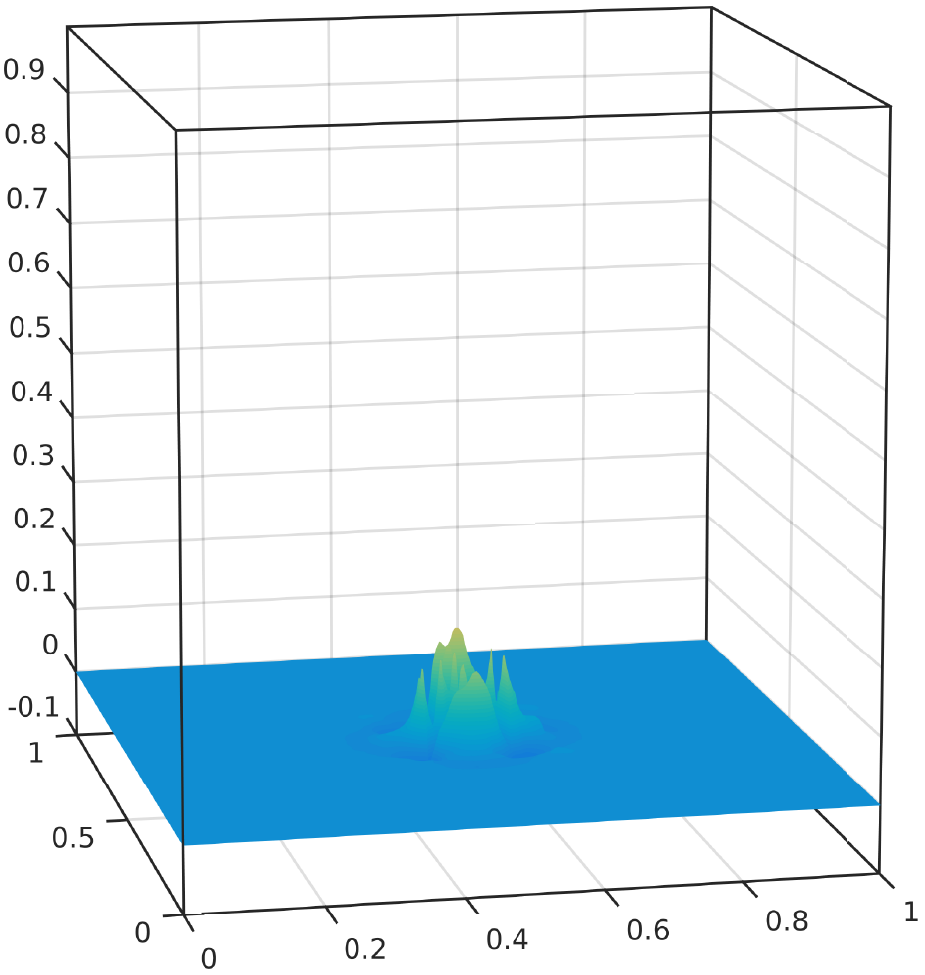}\\
\includegraphics[width=.5\textwidth]{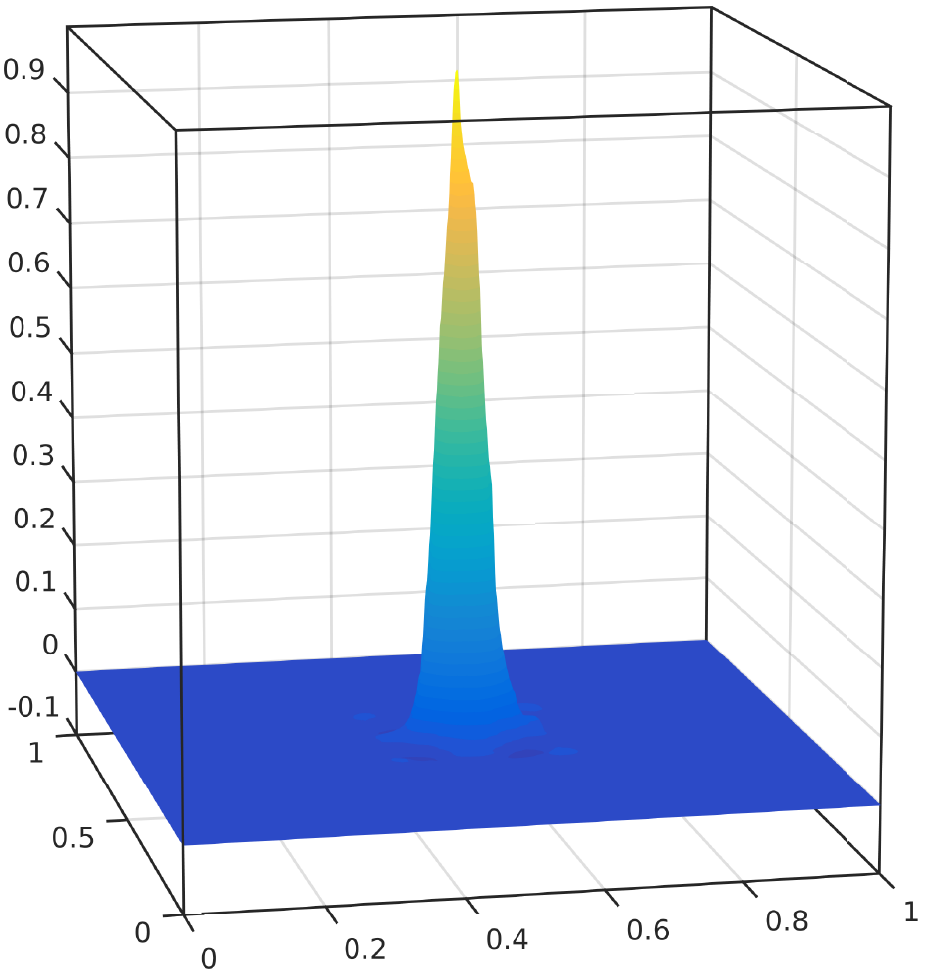}
\includegraphics[width=.5\textwidth]{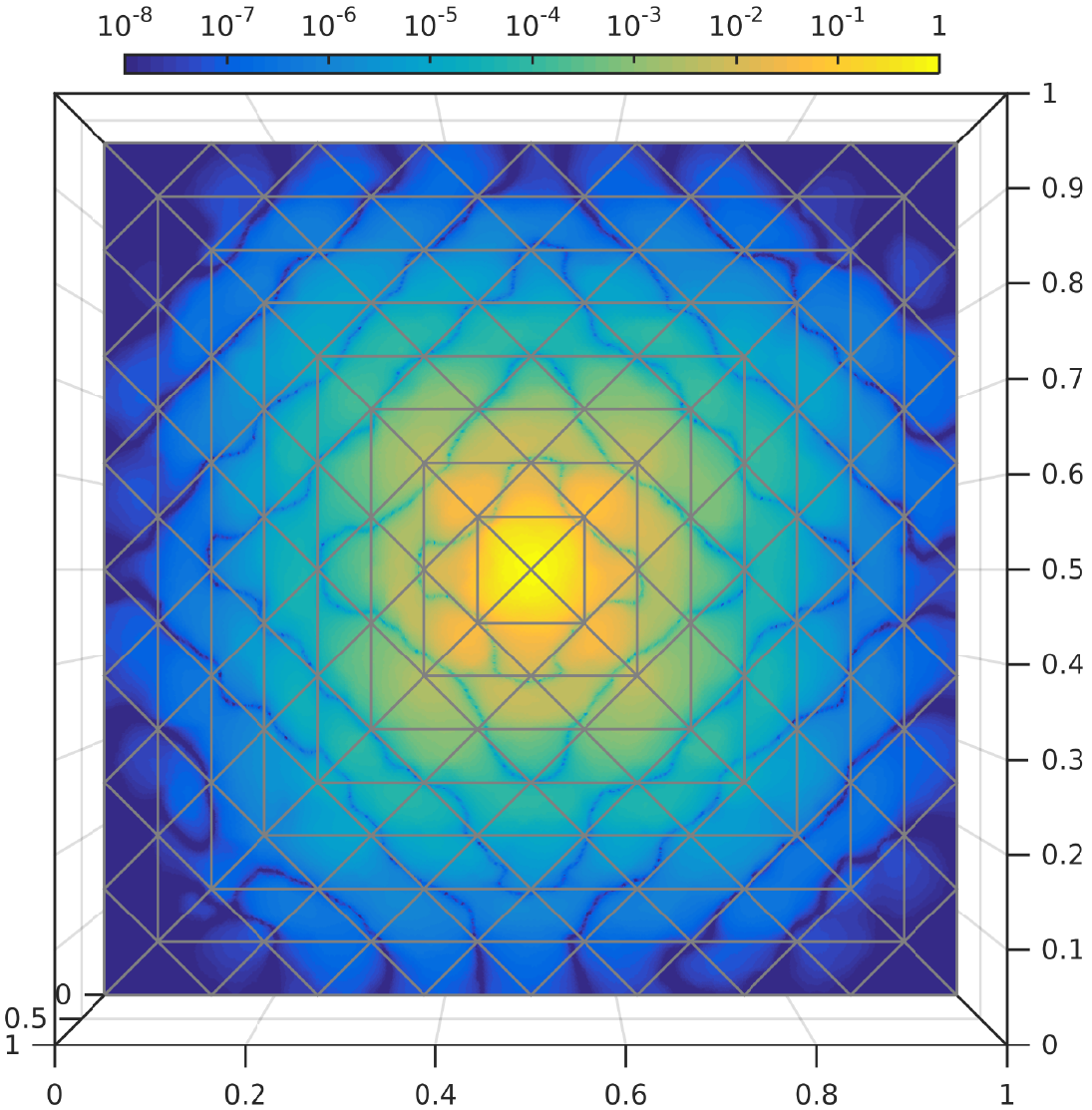}
%\includegraphics[height=.5\textwidth]{LODtest_testl}%
% If no graphics program available, insert a blank space i.e. use
%\picplace{5cm}{2cm} % Give the correct figure height and width in cm
%
\caption{Standard nodal basis function $\lambda_z$ with respect to the coarse mesh $\G_H$ (top left), corresponding ideal corrector $\phi_z=\Cor \lambda_z$ (top right), and corresponding test basis function $\T\lambda_z=(1+\Cor)\lambda_z$ (bottom left). The bottom right figure shows a top view on the modulus of test basis function $\T\lambda_z=(1+\Cor)\lambda_z$ with logarithmic color scale to illustrate the exponential decay property. The underlying rough diffusion coefficient $A$ is depicted in Fig.~\ref{fig:coeff}.}
\label{fig:testbasis_ideal}       % Give a unique label
\end{figure}
We shall now explain how the abstract theory of the previous sections is related to the method of \cite{MP14} and its variants. 
Let $\G_H$ denote some regular (in the sense of Ciarlet) finite element mesh into closed simplices and let $V_H:=P_1(\G_H)\cap V$ denote the space of continuous functions that are affine when restricted to any element $T\in\G_H$. Let $I_H:V\to V_H$ be a quasi-interpolation operator that
acts as a stable quasi-local projection in the sense that
$I_H\circ I_H = I_H$ and that
for any $T\in\G_H$ and all $v\in V$ there holds
\begin{equation}\label{e:IHapproxstab}
H^{-1}\|v-I_H v\|_{L^2(T)} + \|I_H v \|_{V,T}
\leq C_{I_H} \|\nabla v\|_{V,\Omega_T},
\end{equation}
where $\Omega_T$ refers to some neighbourhood of $T$ (typically the union of $T$ and the adjacent elements) and $\|\cdot\|_V:=\|\nabla\cdot\|_{L^2(\Omega)}$. One possible choice (among many others) is to define $I_H:=E_H\circ\Pi_H$, where
$\Pi_H$ is the piecewise $L^2$ projection onto $P^1(\G_H)$ and $E_H$ is the averaging operator that maps $P_1(\G_H)$ to $V_H$ by
assigning to each interior vertex the arithmetic mean of the corresponding
function values of the adjacent elements, that is, for any $v\in P_1(\G_H)$
and any free vertex $z\in\N_H$,
\begin{equation*}
(E_H(v))(z) = \frac{1}{\card\{K\in\G_H\,:\,z\in K\}}\sum_{T\in\G_H:z\in T}v|_T (z).
\end{equation*}
For this choice, the proof of \eqref{e:IHapproxstab} follows from combining the well-established approximation and stability properties of 
$\Pi_H$ and $E_H$, see for example \cite{ern}. The choice of $I_H$ in \cite{MP14,HP13} was slightly different. Therein, the $L^2(\Omega)$-orthogonal projection onto $V_H$  played the role of $I_H$. Since this a non-local operator, the analysis was based on the fact that the local quasi-interpolation operator of \cite[Section 6]{CV} has the same kernel and, hence, induces the same method. 

Following the recipe of Section~\ref{s:PG} and taking into account the present setting with an inner product $a$, the ideal test space
$W_H:=(\kernel I_H)^{\perp_a}$ is simply the orthogonal complement (w.r.t. $a$) of the fine scale functions $\kernel I_H$. 

Given the nodal basis of $V_H$, a basis of $W_H$ is computed by means of the trial-to-test operator $\T=1+\Cor$, where 
\begin{equation}\label{e:correctorhom}
\forall w\in\kernel I_H:\;a(\Cor \lambda_z,w)=-a(\lambda_z,w).
\end{equation}
It is easily checked that the assumptions made in Section~\ref{s:localization} are satisfied in the present setting. In particular, formula \eqref{e:cutclos} holds  with $C_{\eta,I_H}=C_{I_H}(C_{I_H} C_\eta+1)$ and $m=2$. This follows from the product rule, \eqref{e:etaestimate}, and the local approximation and stability properties \eqref{e:IHapproxstab} of $I_H$. This implies the exponential decay as it is stated in \eqref{e:decay0} with constants $C$ and $c$ independent of variations of the diffusion coefficient $A$. An example of a corrector and a test basis function are depicted in Fig.~\ref{fig:testbasis_ideal} to demonstrate the exponential decay.

We truncate the computational domain of the corrector problems to local subdomains of diameter $\ell H$ roughly. We have not yet described how to do this in practice. The obvious way would be to simply replace $\Omega$ in \eqref{e:correctorhom} with suitable neighbourhoods of the nodes $z$. This procedure was used in \cite{MP14}. However, it turned out that it is advantageous to consider the following slightly more involved technique based on element correctors \cite{HP13,HMP14}. 

\begin{figure}[tb]
\sidecaption[t]
\includegraphics[height=0.211\textwidth]{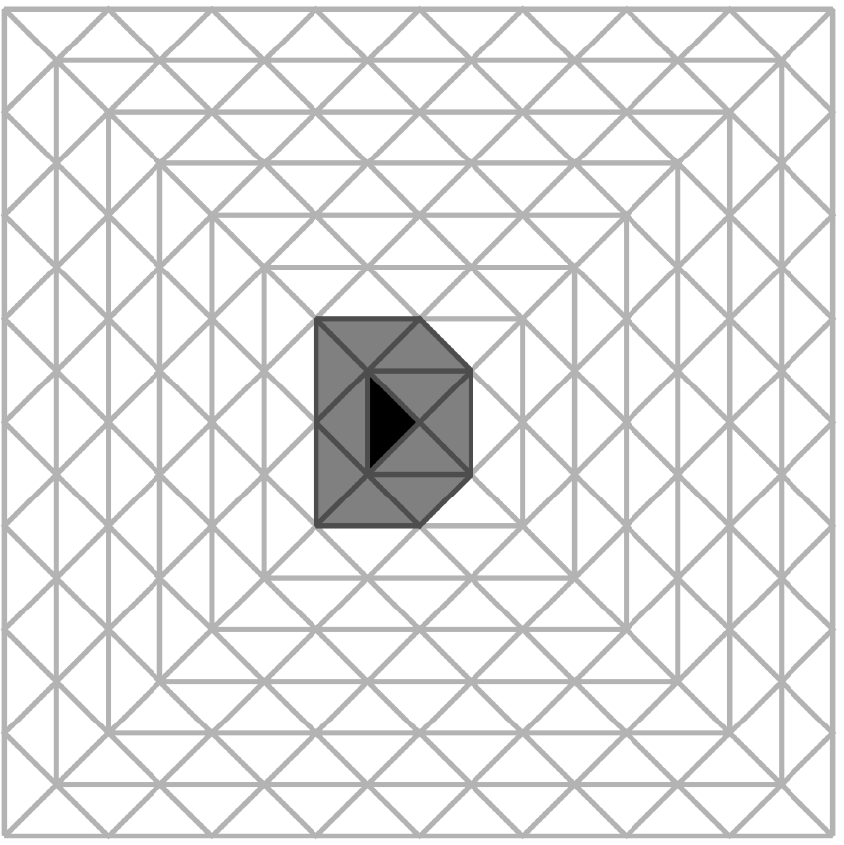}
\includegraphics[height=0.211\textwidth]{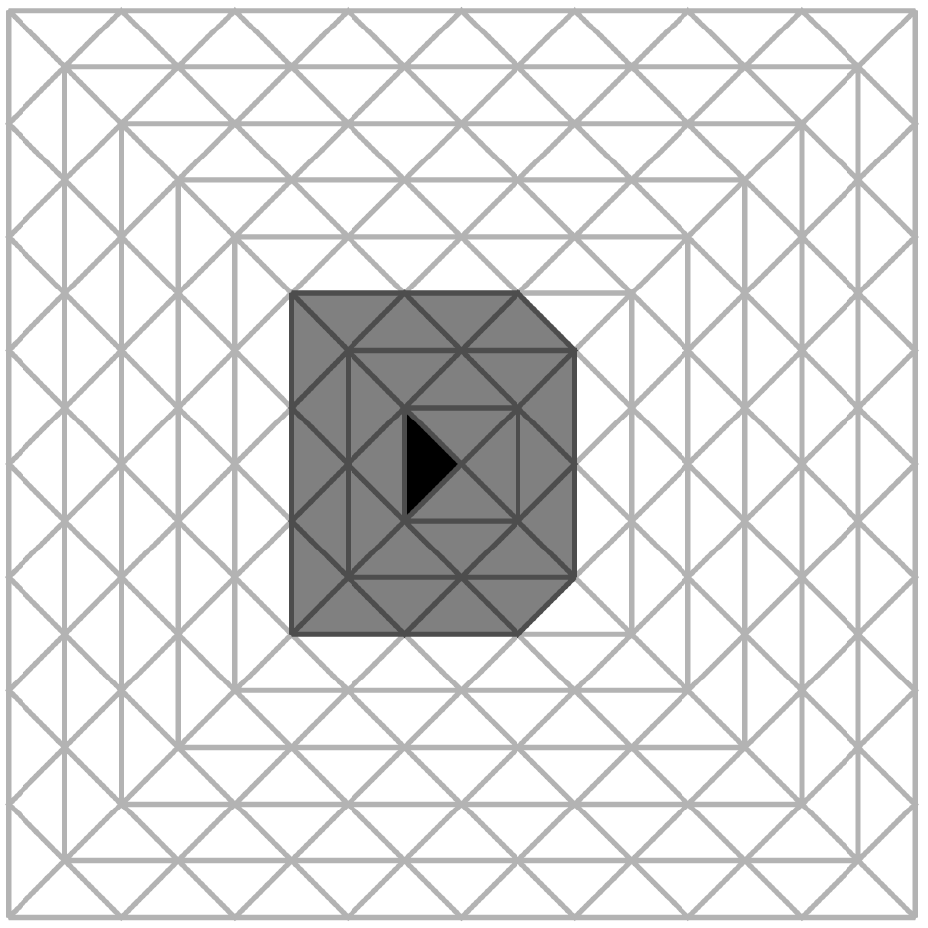}
\includegraphics[height=0.211\textwidth]{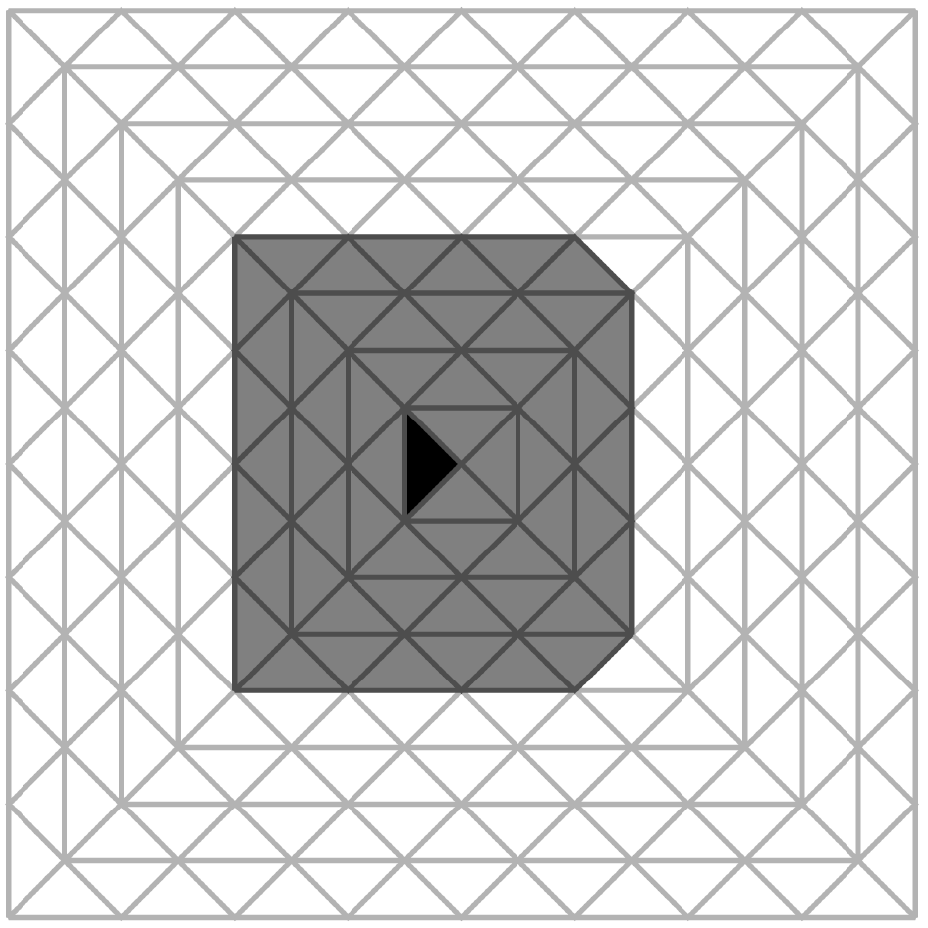}
\caption{Element patches $\Omega_{T,\ell}$ for $\ell=1,2,3$ (from left to right) as they are used in the localised corrector problem \eqref{e:correctorproblemloc}.\label{fig:patches}}
\end{figure}

We assign to any $T\in\G_H$ its $\ell$-th order element patch
$\Omega_{T,\ell}$ for a positive integer $\ell$; see Fig.~\ref{fig:patches} for an illustration. Moreover, we define for all $v,w\in V$ and $\omega\subset\Omega$ the localised bilinear forms
\begin{equation*}
a_{\omega}(v,w):= \int_\omega (A\nabla v)\cdot\nabla w.
\end{equation*}

Given any nodal basis function $\lambda_z\in V_H$,
let $\phi_{z,\ell,T}\in \kernel I_H\cap W^{1,2}_0(\Omega_{T,\ell})$ solve the subscale corrector problem
\begin{equation}\label{e:correctorproblemloc}
a_{\Omega_{T,\ell}}(\phi_{z,\ell,T},w) = -a_T(\lambda_z,w)
\quad\text{for all } w\in \kernel I_H\cap W^{1,2}_0(\Omega_{T,\ell}).
\end{equation}
Let $\phi_{z,\ell}:=\sum_{T\in\G_H:z\in T} \phi_{z,\ell,T}$
and define the test function
\begin{equation*}
\Lambda_{z,\ell} := \lambda_z + \phi_{z,\ell}.
\end{equation*}

\begin{figure}[t]
%\sidecaption[t]
% For example, with the graphicx style use
\includegraphics[width=.245\textwidth]{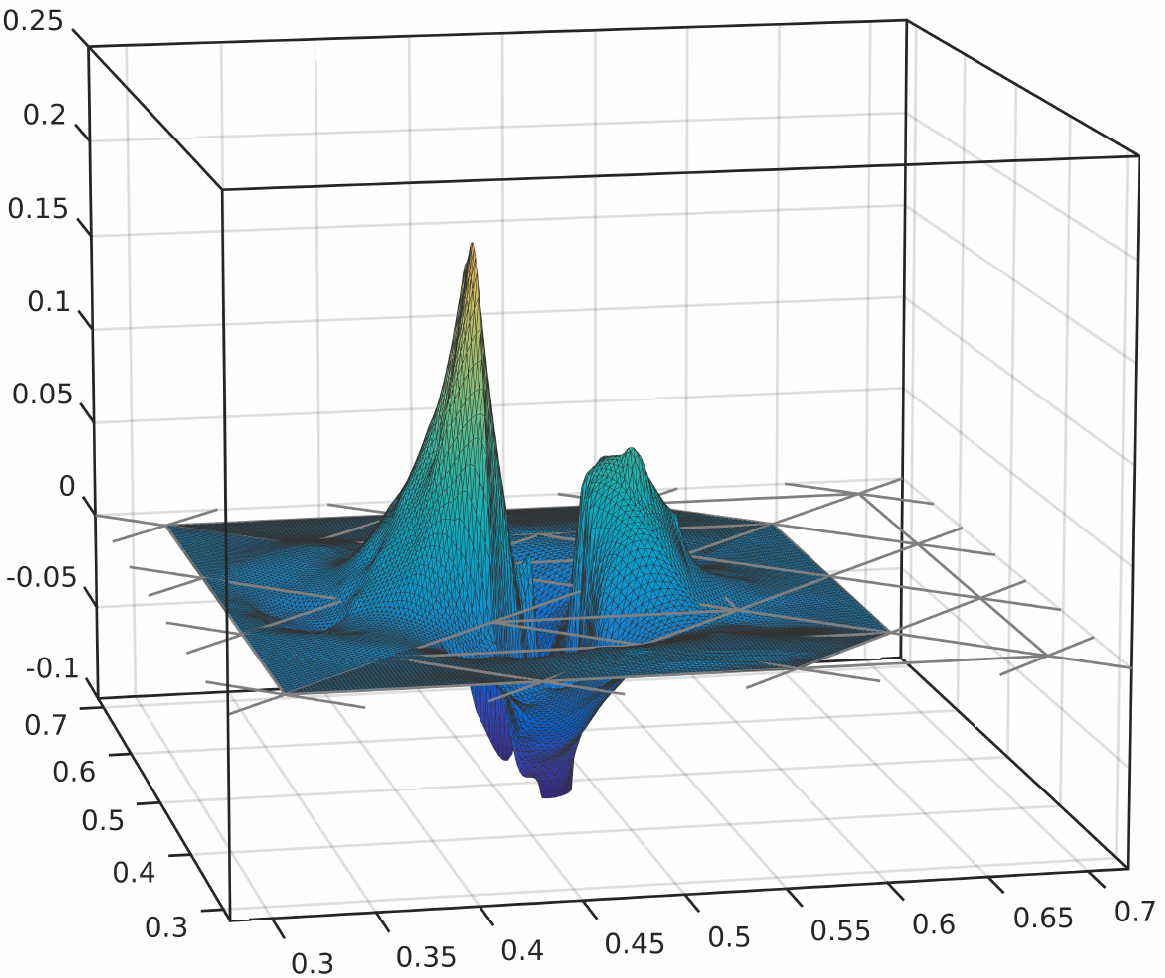}
\includegraphics[width=.245\textwidth]{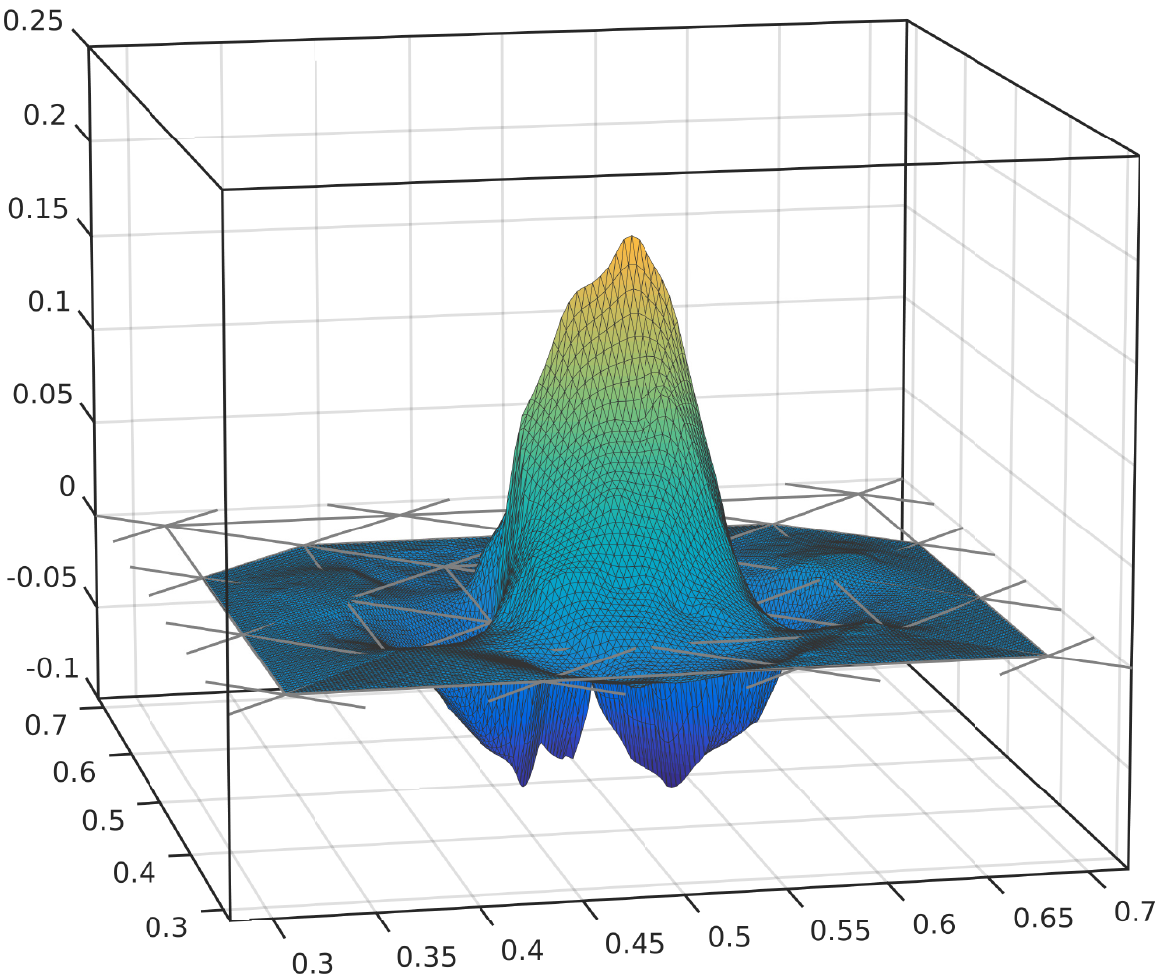}
\includegraphics[width=.245\textwidth]{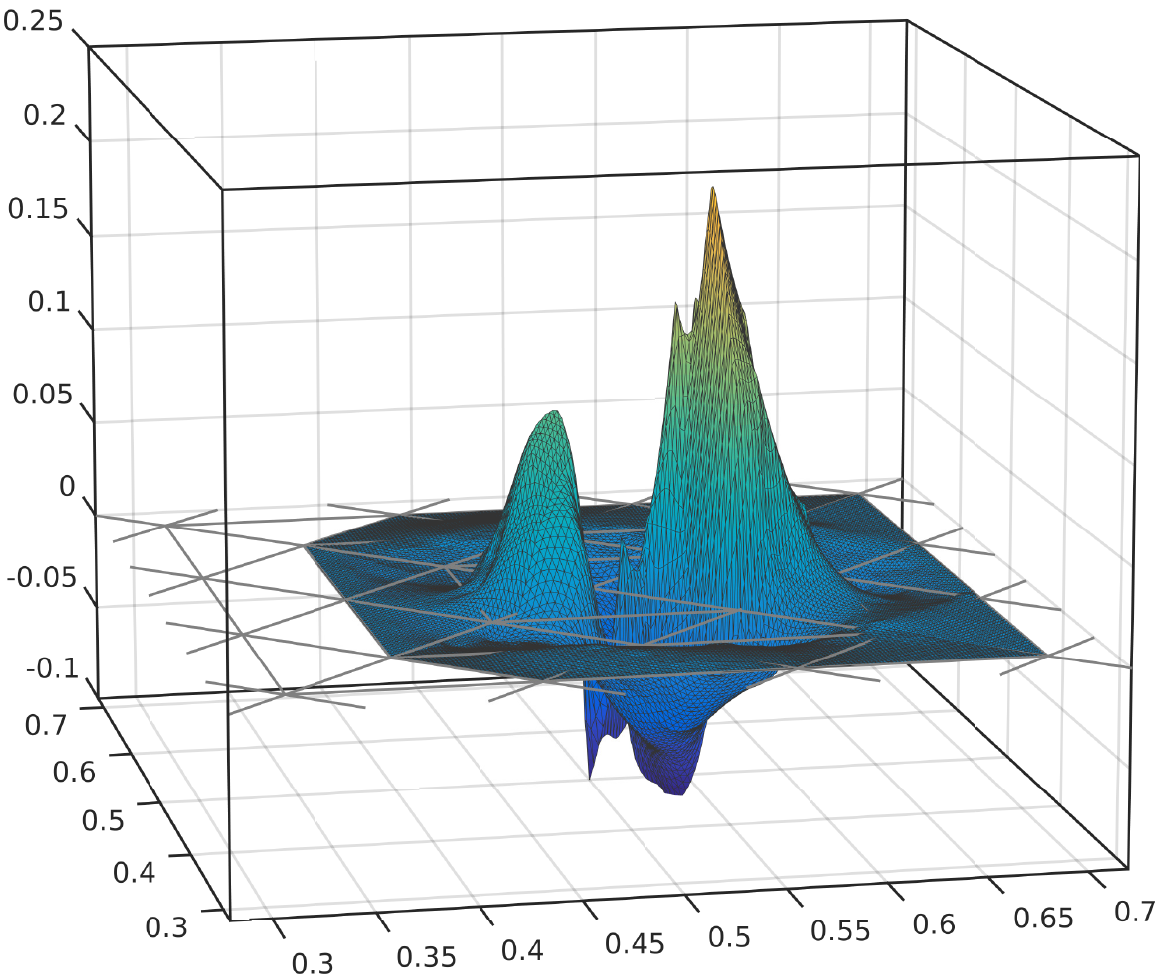}
\includegraphics[width=.245\textwidth]{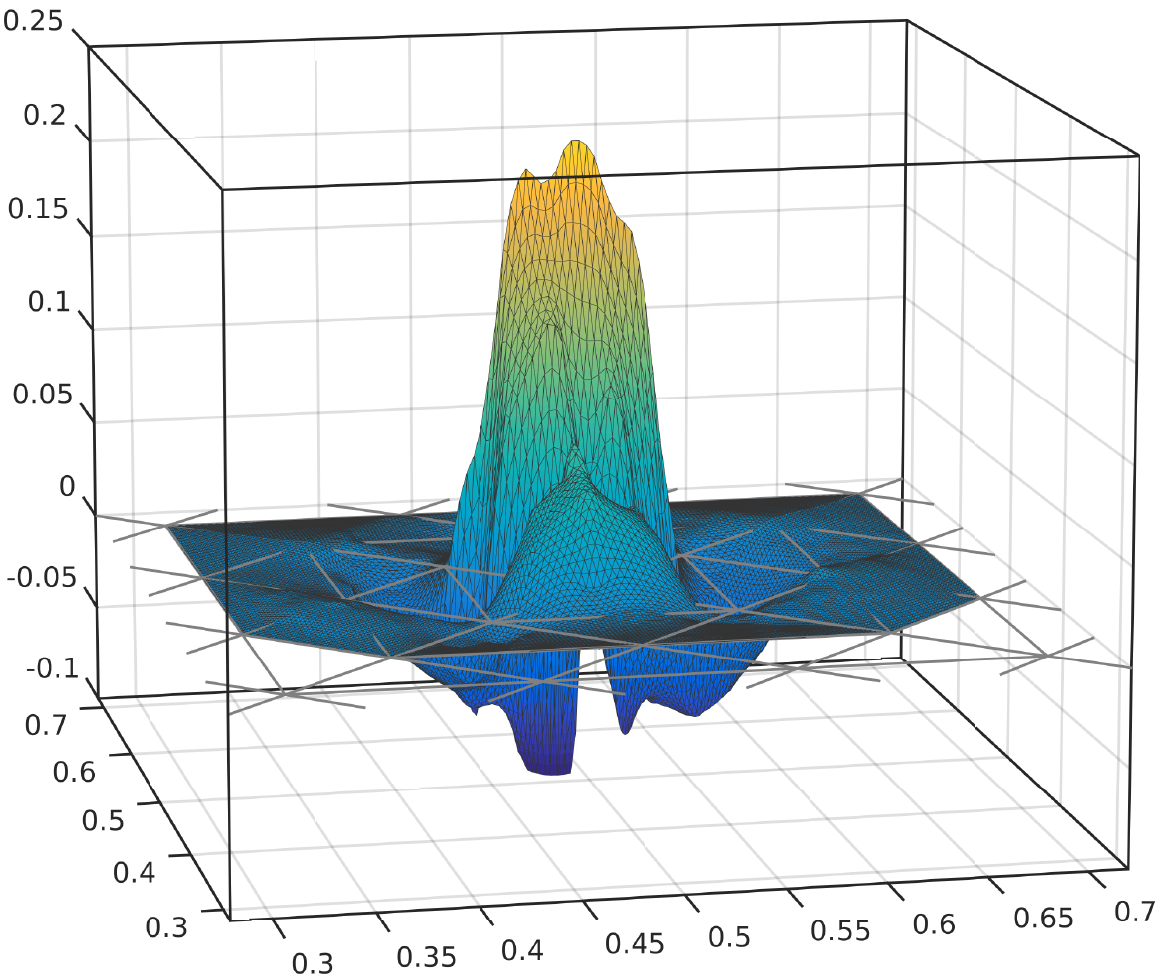}\vspace{1ex}\\
\includegraphics[width=.5\textwidth]{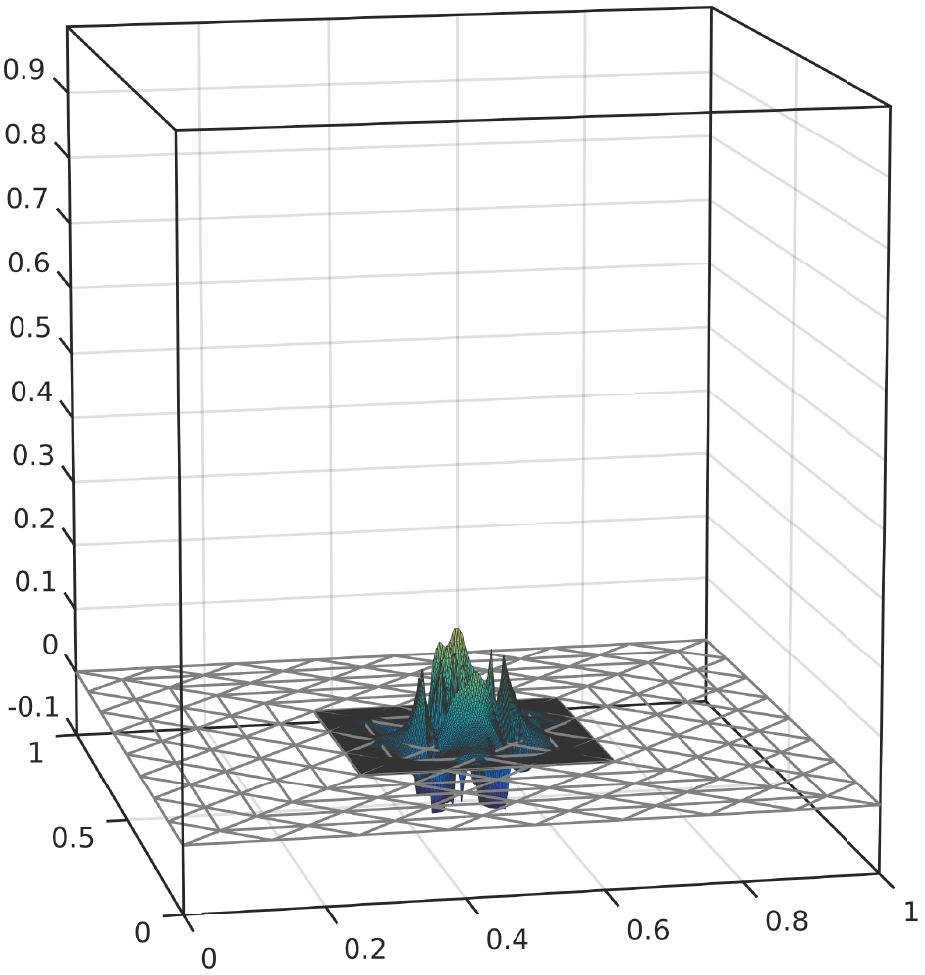}
\includegraphics[width=.5\textwidth]{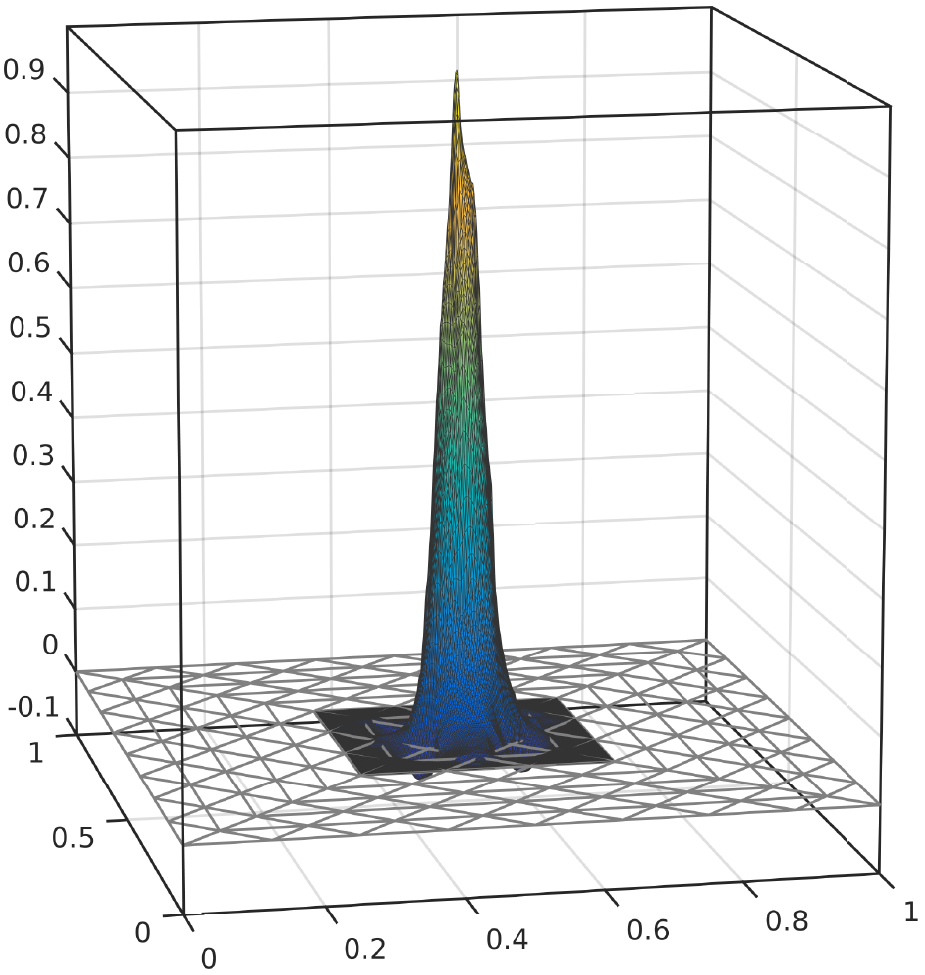}%
% If no graphics program available, insert a blank space i.e. use
%\picplace{5cm}{2cm} % Give the correct figure height and width in cm
%
\caption{Localised element correctors $\phi_{z,\ell,T}$ for $\ell=2$ and all four elements $T$ adjacent to the vertex $z=[0.5,0.5]$ (top), localised nodal corrector $\phi_{z,\ell}=\Cor_\ell\lambda_z=\sum_{T\ni z}\phi_{z,\ell,T}$ (bottom left) and corresponding test basis function $\Lambda_{z,\ell}=\T_{\ell}\lambda_z=(1+\Cor_{\ell})\lambda_z$ (bottom right). The underlying rough diffusion coefficient is depicted in Fig.~\ref{fig:coeff}. The computations have been performed by standard linear finite elements on local fine meshes of with $h=2^{-8}$. See Fig.~\ref{fig:testbasis_ideal} for a comparison with the ideal global  corrector and basis.}
\label{fig:testbasis}       % Give a unique label
\end{figure}

The localised test basis function $\Lambda_{z,\ell}$ and the underlying correctors $\phi_{z,\ell,T}$ can be seen in Fig.~\ref{fig:testbasis}. Note that we impose homogeneous Dirichlet boundary condition on the artificial boundary of the patch which is well justified by the fast decay.

More generally, we may define the localised correction operator $\Cor_\ell$ by 
\begin{equation*}
\Cor_\ell v_H:=\sum_{z\in\N_H}v_H(z)\phi_{z,\ell}
\end{equation*}
as well as the localised trial-to-test operator 
\begin{equation*}
\T_\ell v_H:=1+\Cor_\ell v_H=\sum_{z\in\N_H}v_H(z)\Lambda_{z,\ell}.  
\end{equation*}
The space of test functions then reads
\begin{equation*}
W_H^\ell := \T_\ell V_H = \operatorname{span}\{\Lambda_{z,\ell}\,:\,z\in\N_H\}
\end{equation*}
and the (localised) multiscale Petrov-Galerkin FEM seeks $u_{H,\ell}\in V_H$ such that
\begin{equation}\label{e:discreteproblemhom}
a(u_{H,\ell},w_{H,\ell}) = 
(f,w_{H,\ell})_{L^2(\Omega)}
\quad\text{for all } w_{H,\ell}\in W_{H,\ell}.
\end{equation}
In previous papers \cite{MP14,HP13,HMP14} we have considered the symmetric version with $W_{H,\ell}$ as trial and test space and also the reverse version with $W_{H,\ell}$ as the trial space and $V_H$ as test space \cite{EGH15}. All these methods are essentially equal in the ideal case and there are no major changes in the output after localisation (when only the $V_H$ part of the discrete solution is considered). When it comes to implementation and computational complexity, the present Petrov-Galerkin version has the advantage that there is no communication between the correctors. This means that the fine-scale solutions of the corrector problems need not to be stored but only their interaction with the $\mathcal{O}(\ell^d)$ standard nodal basis functions in their patches; see also \cite{EGH15} for further discussions regarding those technical details. 

The error analysis of the localised method follows similar arguments. Let $u_H\in V_H$ be the ideal Petrov-Galerkin approximation and let $e_H:=u_H-u_{H,\ell}\in V_H$ denote the error with respect to the ideal method.  Then there exists some $z_H\in W_H$ with $\|z_H\|_V=1$ such that
\begin{equation*}
\frac{\alpha}{\CIH}\|e_H\|_V\leq a(e_H,z_H) = a(u_{H,\ell}-u,z_H-z_{H,\ell}),
\end{equation*}
where $z_{H,\ell}\in W_{H,\ell}$. The exponential decay property allows one to choose $z_{H,\ell}$ in such a way that $\|z_H-z_{H,\ell}\|_V\leq \tilde{C}\exp(-c\ell)$; see for instance \cite{HP13,HMP14}. 
This shows that
\begin{eqnarray*}
\|u-u_{H,\ell}\|_V&\leq& \|u-u_H\|_V + \|u_H-u_{H,\ell}\|_V\\
&\leq& \|u-u_H\|_V + \frac{\CIH C_a}{\alpha}\tilde{C}\exp(-c\ell)\|u-u_{H,\ell}\|_V.
\end{eqnarray*}
We shall emphasise that, in the present context, the constants $\tilde{C}$ and $c$ are independent of variations of the rough diffusion tensor but they may depend on the contrast (the ratio between the global upper and lower bound of $A$). Using \eqref{e:quasibestV}, this shows that the moderate choice $\ell\geq|\log(\alpha/(2 \CIH C_a \tilde{C}))|/c=\mathcal{O}(1)$ implies the quasi-optimality (and also the well-posedness) of the Petrov-Galerkin method with respect to the $V$-norm
\begin{equation*}
\|u-u_{H,\ell}\|_V\leq 2\CIH\min_{v_H\in V_H}\|u-v_H\|_V.
\end{equation*}

With regard to the fact that the $V$-best approximation may be poor and standard Galerkin would have provided us with an even better estimate at lower cost, this result is maybe not very impressive. Let us see if we can do something similar for the $L^2$-norm which appears to be the relevant measure in the context of homogenization problems. A standard duality argument shows
that 
\begin{eqnarray*}
\|e_H\|_{L^2(\Omega)}^2&= & a(e_H,z_H)=a(u-u_{H,\ell},z_H-z_{H,\ell})
%C_{I_H}\|e\|_{L^2(\Omega)}\|\Cor e\|_V
\end{eqnarray*}
for some $z_H\in W_H$ with $\|z_H\|_V\leq C_3\alpha^{-1}\CIH\| e_H\|_{L^2(\Omega)}$ and $z_{H,\ell}:=\T_\ell I_H z_H\in W_{H,\ell}$. Similar arguments as before yield
\begin{eqnarray*}
\|u-u_{H,\ell}\|_{L^2(\Omega)}&\leq& C_1\min_{v_H\in V_H}\|u-v_H\|_{L^2(\Omega)}+C_2\exp(-c\ell)\min_{v_H\in V_H}\|u-v_H\|_V,
\end{eqnarray*}
where $C_1:=\|I_H\|_{\mathcal{L}(L^2(\Omega))}$ and $C_2:= C_a\tilde{C}C_3\alpha^{-1}\CIH$.
This shows that the method is accurate also in the $L^2$-norm regardless of the regularity of $u$. If the oversampling parameter is chosen such that $\ell\gtrsim \log H$, then the method is $\mathcal{O}(H)$ accurate in $L^2(\Omega)$ with no pre-asymptotic phenomena. This is the best worst-case rate one can expect for general $f\in V'$ and $A\in L^\infty$.

Note that the previous results hold true for general $L^\infty$-coefficients and all constants are independent of the variations of the diffusion tensor as far as the contrast remains moderately bounded. In particular, the approach is by no means restricted to periodic coefficients or scale separation. For a more detailed discussion of high-contrast problems in this context we refer to \cite{Peterseim.Scheichl:2014}. 
\medskip 
\begin{figure}[t]
\sidecaption[t]
\includegraphics[height=0.3\textwidth]{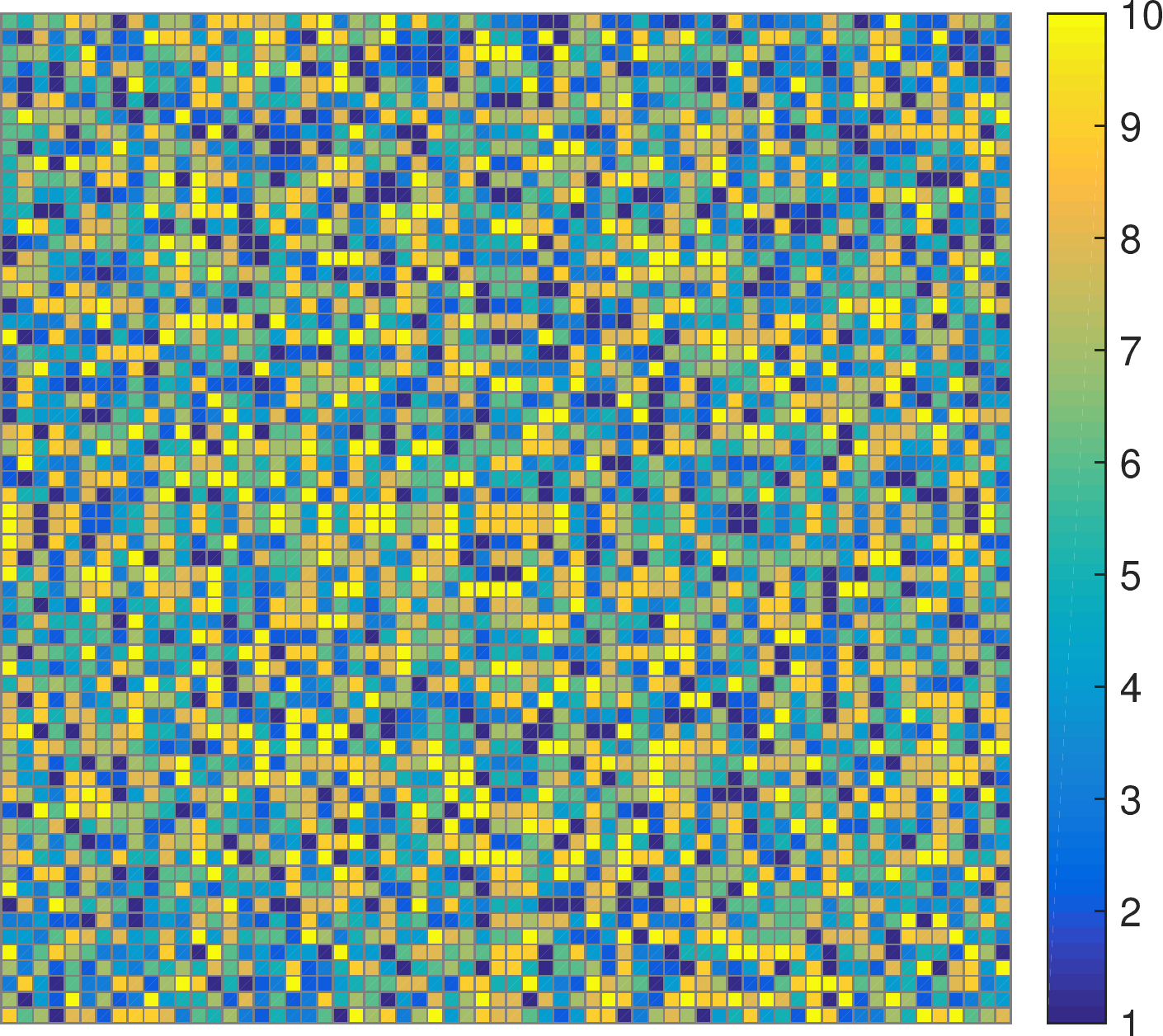}\hspace{1.2ex}
\includegraphics[height=0.295\textwidth]{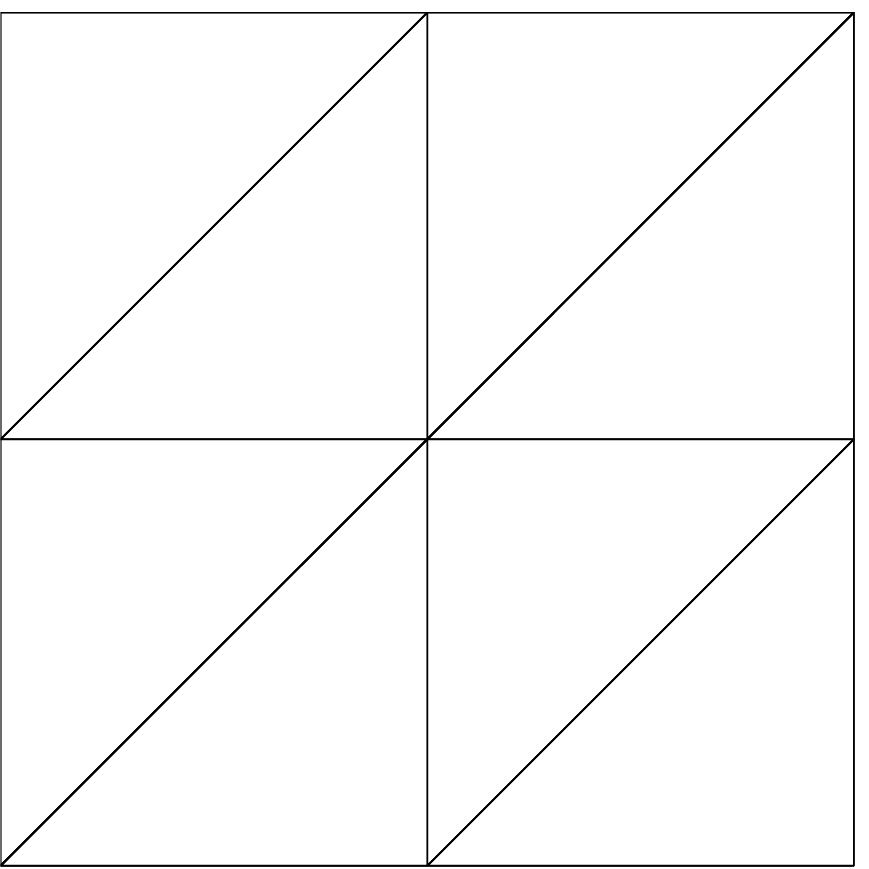}
\caption{Diffusion coefficient in the numerical experiment of Section~\ref{s:hom} and coarsest mesh.\label{fig:coeff}}
\end{figure}
\begin{figure}[t]
\sidecaption[t]
% For example, with the graphicx style use
\includegraphics[width=.6\textwidth]{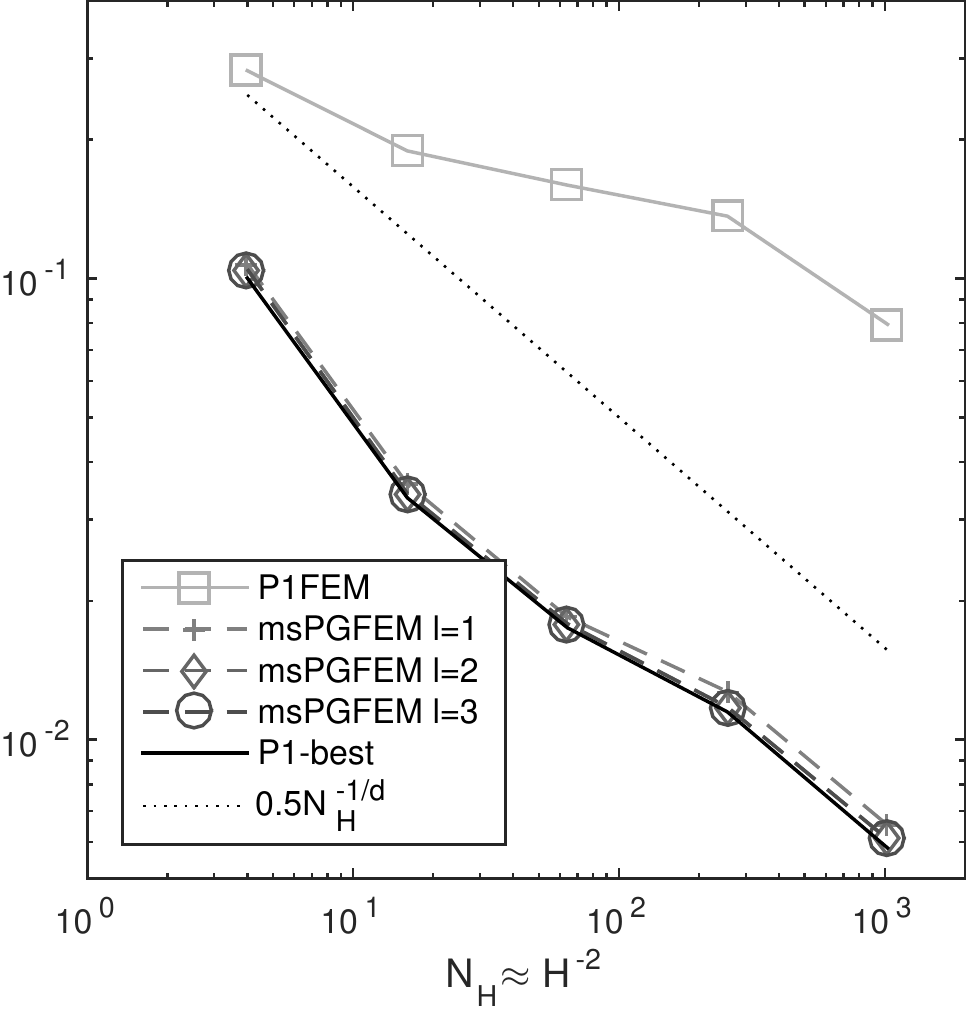}
%
% If no graphics program available, insert a blank space i.e. use
%\picplace{5cm}{2cm} % Give the correct figure height and width in cm
%
\caption{Numerical experiment of Section~\ref{s:hom}. Relative $L^2$-errors of multiscale Petrov-Galerkin FEM \eqref{e:discreteproblemhom} versus the number of degrees of freedom $N_H\approx H^{-2}$, where $H= 2^{-1},\ldots,2^{-5}$ is the uniform coarse mesh size. The localisation parameter varies between $\ell=1,\ldots,3$. The $P_1$-FE solution and the best-approximation in the $P_1$-FE space on the same coarse meshes are depicted for comparison.}
\label{fig:expnumhom}       % Give a unique label
\end{figure}

The final step towards a fully practical method is the discretization of the fine-scale corrector problems. With regard to the possible low regularity of the solution, $P_1$ finite elements on a refined mesh $\G_h$ appears reasonable, but any other type of discretisation is possible. Obviously, the fine-scale discretization parameter $h$ has to be chosen fine enough to resolve all relevant features of the diffusion coefficient. The previous theory can be transferred to this case in a straight-forward way and we refer to \cite{MP14,HP13,HM14} for the technical details. 
\medskip

%To illustrate the previous estimates, consider the case of an oscillatory but smooth scalar coefficient $A_\varepsilon$, where $\varepsilon$ represents the characteristic period length such that $\|\nabla A\|_{L^\infty(\Omega)}\leq C\varepsilon^{-1}$. Then the error of the multiscale Petrov-Galerkin method satisfies  
%\begin{equation*}
% \|u-u_{H,\ell,h}\|_{L^2(\Omega)}\leq C_{f,\alpha}\left( \frac{H^2}{\varepsilon} + \frac{h}{\varepsilon}\right),
%\end{equation*}
%whereas the error of the standard FEM will be of order one as long as $H\gtrsim \varepsilon$.

To illustrate the previous estimates, we close this section with a numerical experiment. Let $\Omega$ be the unit square and the outer force $f\equiv 1$ in $\Omega$. Consider the coefficient $A$ that is piecewise constant with respect to a uniform Cartesian grid of width $2^{-6}$. Its values are randomly chosen between $1$ and $10$; see Fig.~\ref{fig:coeff}. Consider uniform coarse meshes $\G_H$ of size $H=2^{-1},2^{-2},\ldots,2^{-5}$ of $\Omega$ that certainly do not resolve the rough coefficient $A$ appropriately. The reference mesh $\G_h$ has width $h=2^{-9}$. Since no analytical solutions are available, the standard finite element approximation $u_h\in V_h$ on the reference mesh $\G_h$ serves as the reference solution. Doing this, we assume that $u_h$ is sufficiently accurate and, necessarily, that $\G_h$ resolves the discontinuities of $A$. The corrector problems are also are also solved on this scale of numerical resolution.

The numerical results, i.e. errors with respect to the reference solution $u_h$ are depicted in Fig.~\ref{fig:expnumhom}. The results are in agreement with the theoretical results. They are even better in the sense that $\ell=1$ seems to be sufficient for quasi-optimality (with respect to $u_h$) in the present setup and parameter regime. We expect that the true errors with respect to $u$ would behave similar in the beginning but level off at some point when the reference error starts to dominate the upscaling error. Still, the experiment clearly indicates that numerical homogenization is possible for very general $L^\infty$-coefficients. 

We refer to \cite{MP14,HP13,HMP14,HMP13,EGH15,AH14,Brown.Peterseim:2014,EGMP13,HM14,MP12,Henning.Mlqvist.Peterseim:2013} for many more numerical experiments for several model problems including nonlinear stationary and non-stationary problems as well as eigenvalue problems. 

\section{Application to high-frequency acoustic scattering}\label{s:helmholtz}
This section will show that the abstract framework of Section~\ref{s:abstract}--\ref{s:localization} is indeed applicable beyond the coercive and symmetric model problem of the previous section. We consider the scattering of acoustic waves at a sound-soft scatterer modelled by the Helmholtz equation over a bounded Lipschitz domain $\Omega\subset\R^{d}$ ($d=1,2,3$),\begin{subequations}
\label{e:modelhelm}
\end{subequations}
\begin{equation}\label{e:modela}
  -\Delta u - \k^2 u = f\quad \text{in }\Omega,
 \tag{\ref{e:modelhelm}.a}
\end{equation}
along with mixed boundary conditions of Dirichlet and Robin type
\begin{align}
  u &= 0\quad\text{on }\Gamma_D,\tag{\ref{e:modelhelm}.b}\label{e:modelb}\\
  \nabla u\cdot \nu - i\k u &= 0\quad\text{on }\Gamma_R.\tag{\ref{e:modelhelm}.c}\label{e:modeld}
\end{align}
Here, the wave number $\k\gg 1$ is real and positive, $i$ denotes the imaginary unit and $f\in L^2(\Omega,\mathbb{C})$. We assume that the boundary $\Gamma :=\partial \Omega $
consists of two components 
\begin{equation*}
\partial\Omega =\overline{\Gamma _D\cup \Gamma _R},\quad\overline{\Gamma}_D\cap\overline{\Gamma}_R =\emptyset
\end{equation*}
where $\Gamma _D$ encloses the scatterer and $\Gamma _R$ is an artificial truncation of the whole unbounded space. The vector $\nu$ denotes the unit normal vector that is outgoing from $\Omega$.

Given $f\in L^2(\Omega,\mathbb{C})$, we wish to find $u\in V:=\{v\in W^{1,2}(\Omega,\C)\;\vert\; v=0\text{ on }\Gamma_D\}$ such that, for all $v\in V$, 
\begin{equation}\label{e:modelvarhelm}
 a(u,v) := \int_\Omega\nabla u\cdot\nabla \bar v-\k^2\int_\Omega u \bar v - i\k\int_{\Gamma_R}u\bar v = \int_\Omega f \bar v=:\overline{F(w)}.
\end{equation}
The space $V$ is equipped with the usual $\k$-weighted norm $\|v\|_V^2:=\k^2\|v\|_{L^2(\Omega)}^2+\|\nabla v\|_{L^2(\Omega)}^2$. The presence of the Robin boundary condition \eqref{e:modeld} ensures that this variational problem is well-posed in the sense of \eqref{e:BNB1} with $\alpha = 1/\Cstab(\k)$ for some $\k$-dependent stability constant $\Cstab(\k)$; see e.g. \cite{MelenkEsterhazy}. The dependence on the wave number $\k$ is not known in general. An exponential growth  with respect to the wave number is possible \cite{betcke} in non-generic domains. In most cases, the growth seems to be only polynomially, although this is an empirical rather than a theoretical statement, and sufficient geometric conditions for this to hold are rare \cite{MelenkEsterhazy,melenk_phd,feng,makridakis}. For the above scattering problem, we know that $\Cstab(\k)\leq\mathcal{O}(\k)$ if $\Omega$ is convex and if the scatterer is star-shaped \cite{hetmaniuk}.

It is this $\k$-dependence in the stability of the problem that makes the numerical approximation by FEM or related schemes extremely difficult in the regime of large wave numbers. Any perturbation of the problem, e.g. by some discretization, can be amplified by $\Cstab(\k)$. We have seen in the introduction that this is indeed observed in practice and causes a pre-asymptotic effect known as the pollution effect or numerical dispersion 
\cite{BabSau}. This effect puts very restrictive assumptions on the smallness of the underlying mesh that is much stronger than the minimal requirement for a meaningful representation of highly oscillatory functions from approximation theory, that is, to have at least $5-10$ degrees of freedom per wave length and coordinate direction. 

It is the aim of many newly developed methods to overcome or at least to reduce the pollution effect; see \cite{MR2219901,MR2551150,MR2813347,perugia,Zitelli20112406,dpg} among many others. However, the only theoretical results regard high-order FEMs with the polynomial degree $p$ coupled to the wave number $\k$ via the relation $p\approx\log\k$ \cite{MS10,MS11,parsania,MelenkEsterhazy}. Under this moderate assumption, those methods are stable and quasi-optimal in the regime $H\k/p\lesssim 1$ for certain model Helmholtz problems. 

The multiscale method of \cite{Peterseim2014} then showed that pollution in the numerical approximation of the Helmholtz problem can also be cured for a fairly large class of Helmholtz problems, including the acoustic scattering from convex non-smooth objects, by stabilization in the present framework. If the data of the problem (domain, boundary condition, force term) allows for polynomial-in-$\k$ bounds of $\Cstab(\k)$ and if the resolution condition $H\k\lesssim 1$ and the oversampling condition $\log(\k)/\ell\lesssim 1$ are satisfied, then the method is stable and quasi-optimal in the $V$-norm. 

The recent paper \cite{Gallistl.Peterseim:2015} interprets the method of \cite{Peterseim2014} in the present framework and we recall it here very briefly. Given the same discrete setup as in the previous section with some simplicial mesh $\G_H$, corresponding $P_1$ FE space $V_H:=P_1(\G_H)\cap V$, and quasi-interpolation operator $I_H:V\rightarrow V_H$, the multiscale Petrov-Galerkin method is formally defined in the same way. We simply replace the inner product of Section~\ref{s:hom} with the sesquilinear form $a$ of this section.  

\begin{figure}[tb]
%\sidecaption[t]
% For example, with the graphicx style use
\includegraphics[width=.33\textwidth]{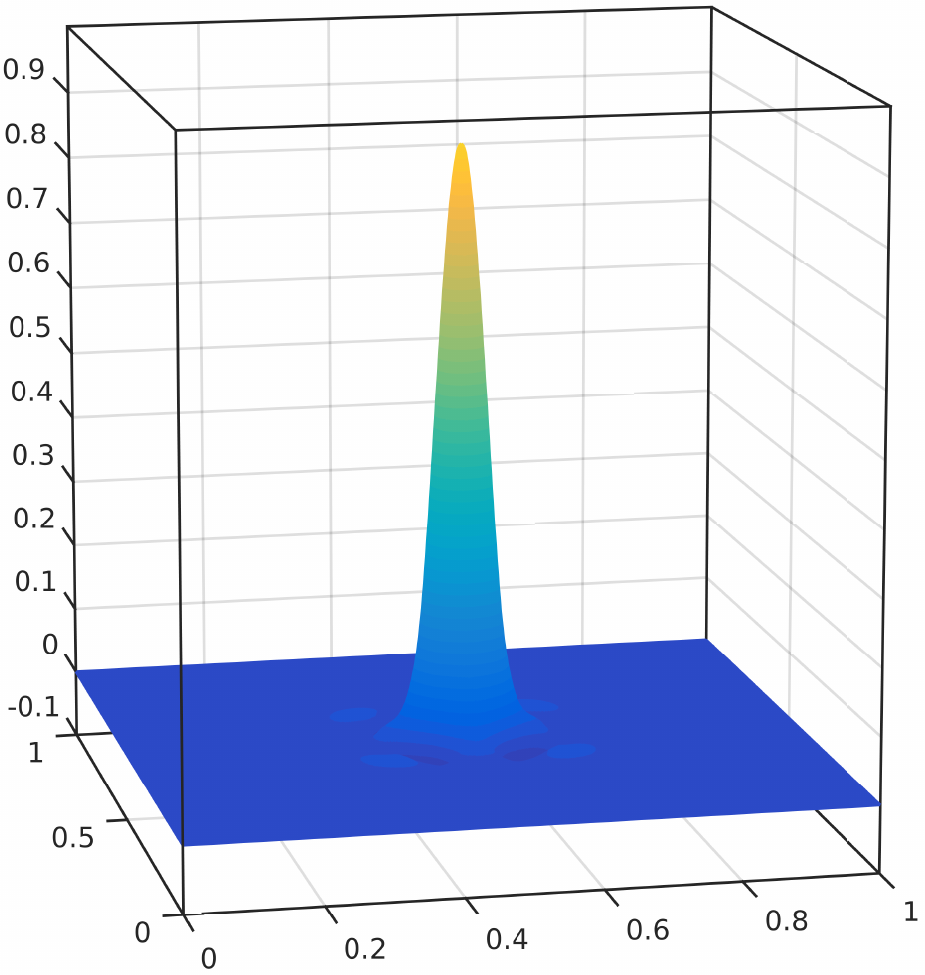}
\includegraphics[width=.33\textwidth]{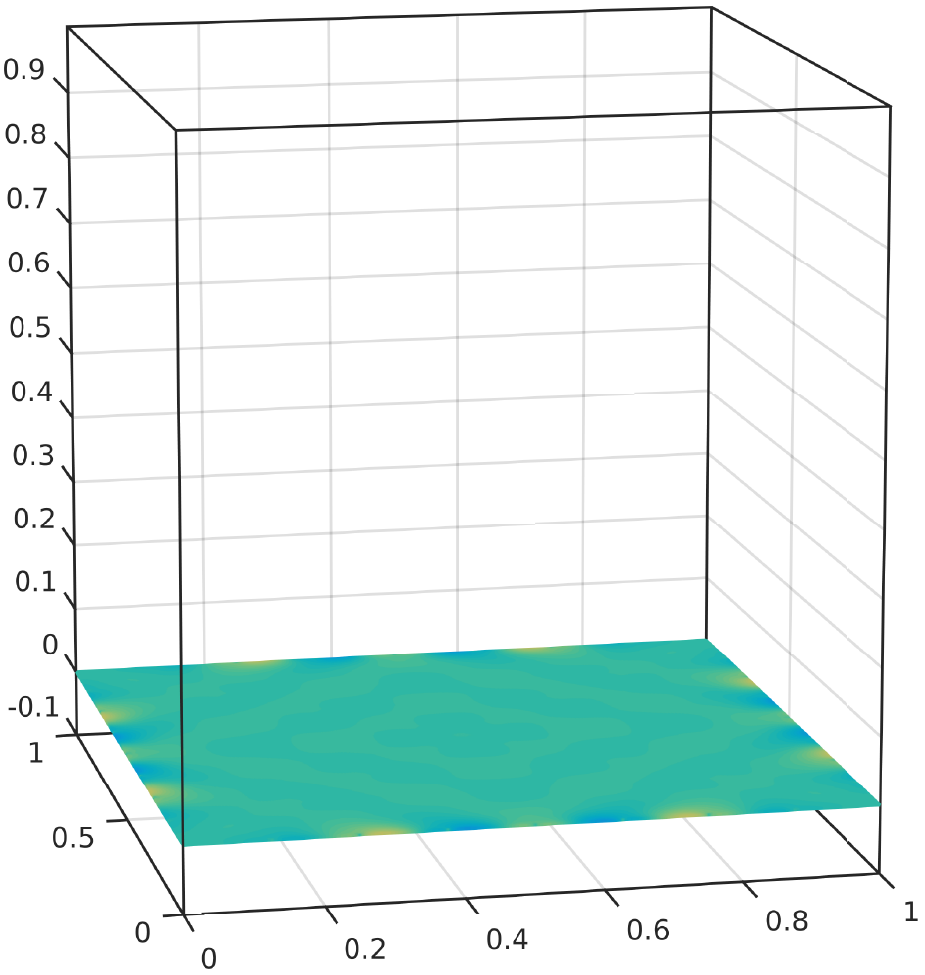}
\includegraphics[width=.33\textwidth]{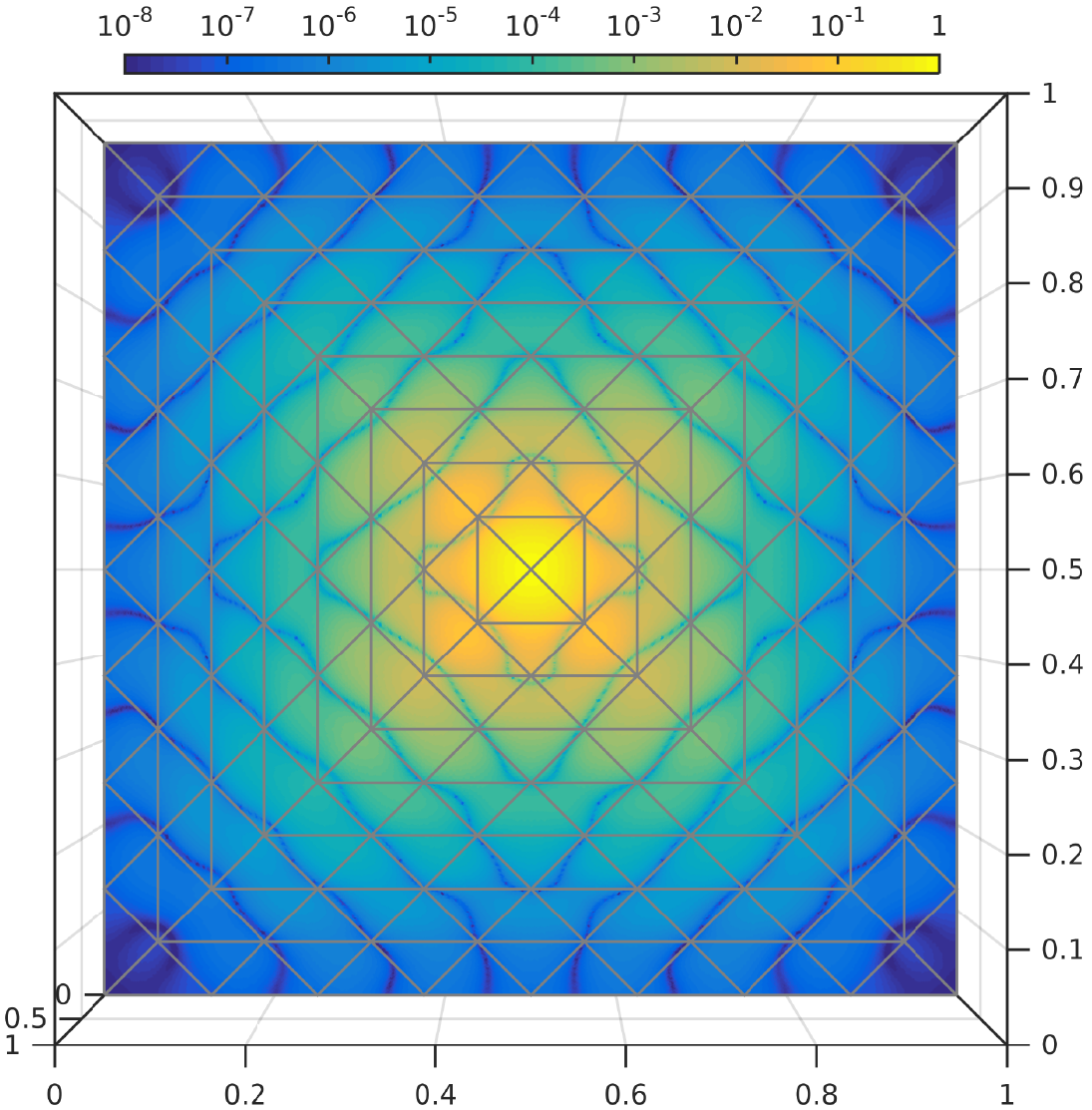}\vspace{1ex}\\
\includegraphics[width=.5\textwidth]{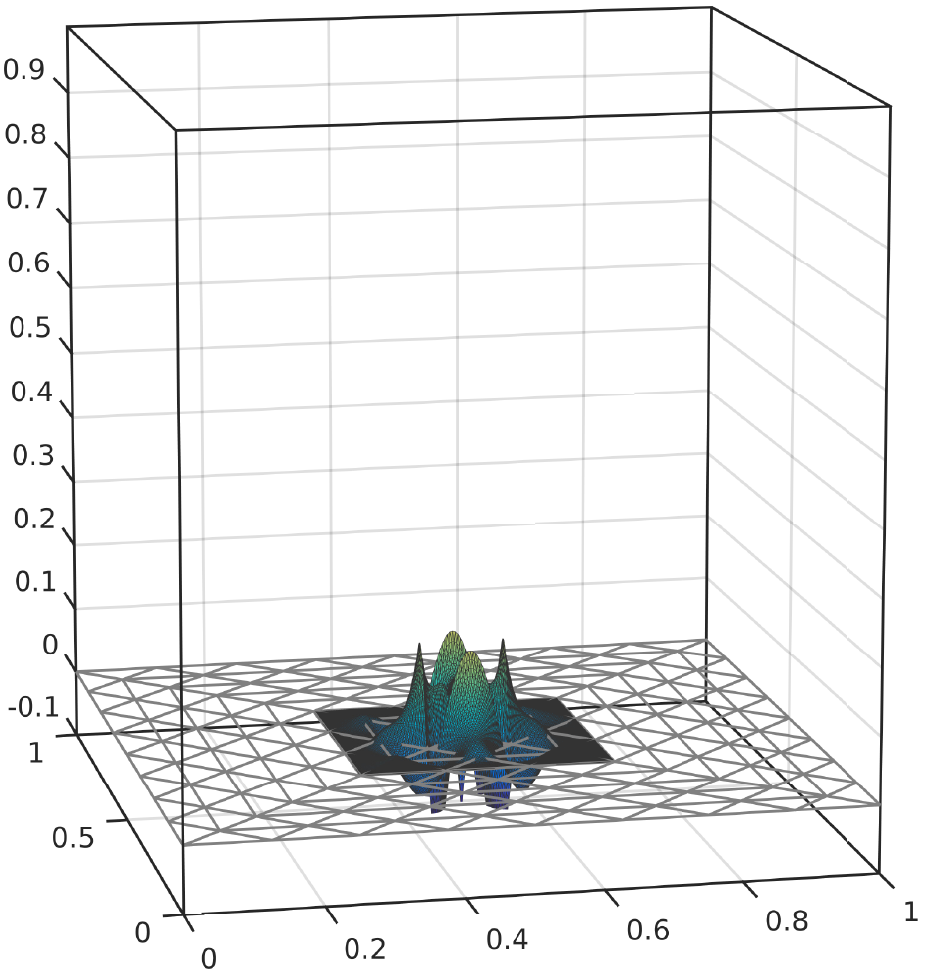}
\includegraphics[width=.5\textwidth]{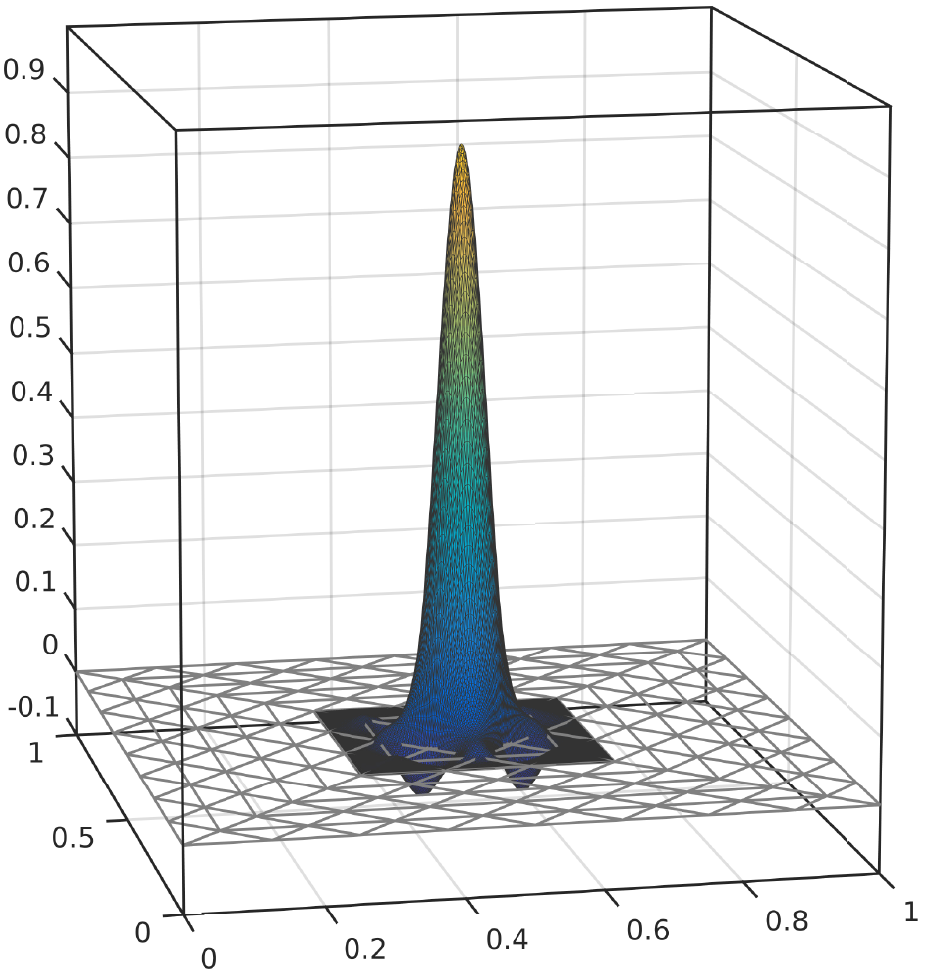}%
% If no graphics program available, insert a blank space i.e. use
%\picplace{5cm}{2cm} % Give the correct figure height and width in cm
%
\caption{Real and imaginary part of the ideal corrector $\Cor \lambda_z$ (top left and middle). The top right figure shows a top view on the modulus of test basis function $\T\lambda_z=(1+\Cor)\lambda_z$ with logarithmic color scale to illustrate the exponential decay property. The underlying computational domain is the unit square with a Robin boundary condition everywhere. The wave number $\k=2^4$ is chosen such that the resolution condition on the coarse mesh is just satisfied. The localised nodal corrector $\phi_{z,\ell}=\Cor_\ell\lambda_z$ (bottom left) and corresponding test basis function $\Lambda_{z,\ell}=\T\lambda_z$ (bottom right) are real-valued because the patch boundary doesn't touch the domain boundary. The local fine meshes used in the computation have width $h=2^{-8}$.}
\label{fig:testbasisHelmholtz}       % Give a unique label
\end{figure}

Given any nodal basis function $\lambda_z\in V_H$, we  construct a corresponding test basis function $\Lambda_{z,\ell}$ by the same procedure as in the previous section, $\Lambda_{z,\ell} := \lambda_z + \phi_{z,\ell}$, where $\phi_{z,\ell}:=\sum_{T\in\G_H:z\in T}\phi_{z,T}$ and $\phi_{z,T}$ solves the cell problem
\begin{equation*}%\label{e:lambdacorrectorproblemhelm}
a_{\Omega_{T,\ell}}(w,\phi_{z,T}) = -a_T(w,\lambda_z)
\quad\text{for all } w\in \kernel I_H\text{ with }\support w\subset\bar{\Omega}_T.
\end{equation*}
Here, 
\begin{equation*}
a_\omega(u,v) := \int_{\Omega\cap\omega}\nabla u\cdot\nabla \bar v-\k^2\int_{\Omega\cap\omega} u \bar v - i\k\int_{\Gamma_R\cap\partial\omega}u\bar v
\end{equation*}
for $\omega \in\{\Omega_{T,\ell},T\}$. Note that the corrector problem inherits the boundary condition from the original problem when the patch boundary coincides with the boundary of $\Omega$. On the part of the patch boundary that falls in the interior of $\Omega$, we simply put the homogeneous Dirichlet condition. 
A major observation is that this corrector problem is well-posed and, in particular, coercive with $\beta=1/3$ under the condition $H\k \leq c_{\operatorname{res}}$ for some given constant $0<c_{\operatorname{res}}=\mathcal{O}(1)$ that only depends on the constant in \eqref{e:IHapproxstab} but not on $H$ or $\k$. This is because $a$ satisfies a G{\aa}rding inequality and fine-scale functions satisfy $\|w\|_{L^2(\Omega)}\leq C_{I_H}H\|\nabla w\|_{L^2(\Omega)}$. 
This coercivity also implies the desired exponential decay of the ideal correctors so that the choice $\Omega_{T,\ell}$ is well justified. This can also be observed in Fig.~\ref{fig:testbasisHelmholtz}. 

The space of localised test functions then reads $
W_{H,\ell} := \operatorname{span}\{\Lambda_{z,\ell}\,:\,z\in\N_H\}$ and the multiscale Petrov-Galerkin FEM seeks $u_{H,\ell}\in V_H$ such that
\begin{equation}\label{e:discreteproblemhelm}
a(u_{H,\ell},w_{H,\ell}) = 
\overline{F(w_{H,\ell})}
\quad\text{for all } w_{H,\ell}\in W_{H,\ell}.
\end{equation}

The quasi-optimality result of the previous section is easily transferred to the present setting. The resolution condition $H\kappa \leq c_{\operatorname{res}}$ and the oversampling condition 
\begin{equation*}
\ell\geq |\log(\alpha/(2 \CIH C_a \tilde{C}))|/c=\mathcal{O}(\log\Cstab(\k))
\end{equation*}
imply the quasi-optimality (and stability) of the multiscale Petrov-Galerkin method
with respect to the $V$-norm
\begin{equation*}
\|u-u_{H,\ell}\|_V\leq 2\CIH\min_{v_H\in V_H}\|u-v_H\|_V.
\end{equation*}
Here, the constants $c$ and $\tilde{C}$ are related to the exponential decay of the test basis (cf. \eqref{e:decay0}) and they are independent of $\k$ under the resolution condition. We shall emphasise that such a best-approximation property does not hold for standard FEMs which require e.g. $\kappa^2H\lesssim 1$ for quasi-optimality \cite{melenk_phd} in the case of pure Robin boundary conditions on a convex planar domain. The FEM approximation is not even known to exist unless $\kappa^{3/2}H\lesssim 1$ in the simplest model problem without a scatterer \cite{Wu2014CIPFEM}.
\begin{figure}[tb]
\sidecaption[t]
\includegraphics[height=0.3\textwidth]{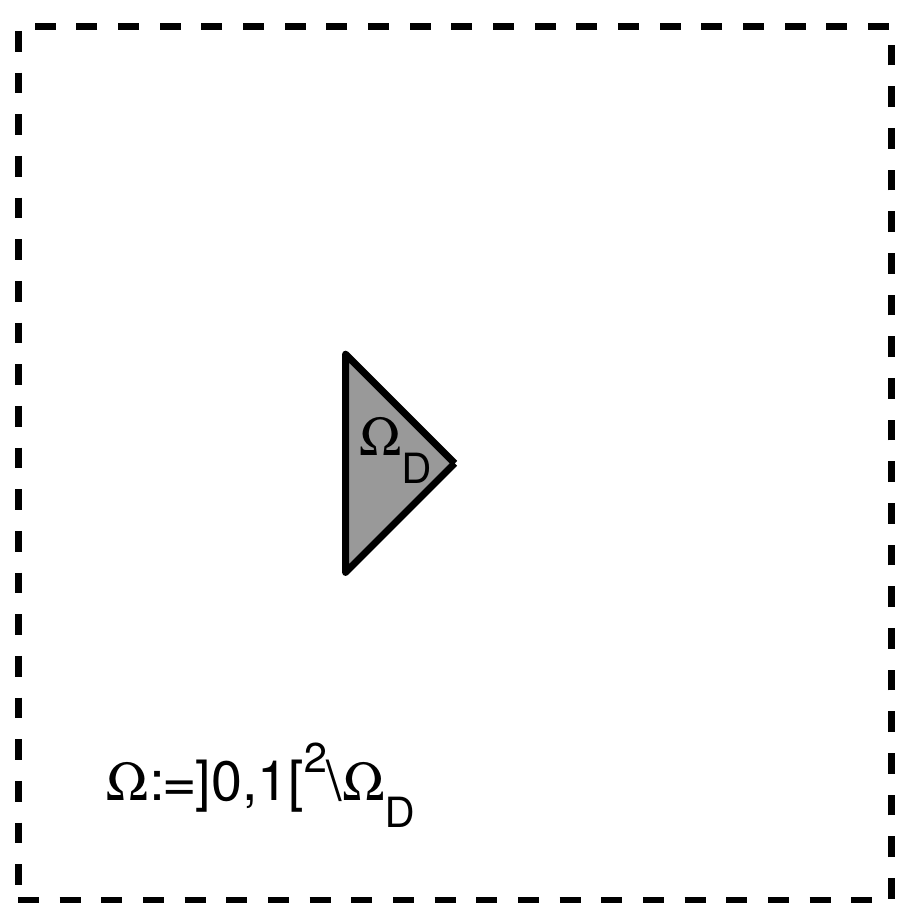}\hspace{1.2ex}
\includegraphics[height=0.3\textwidth]{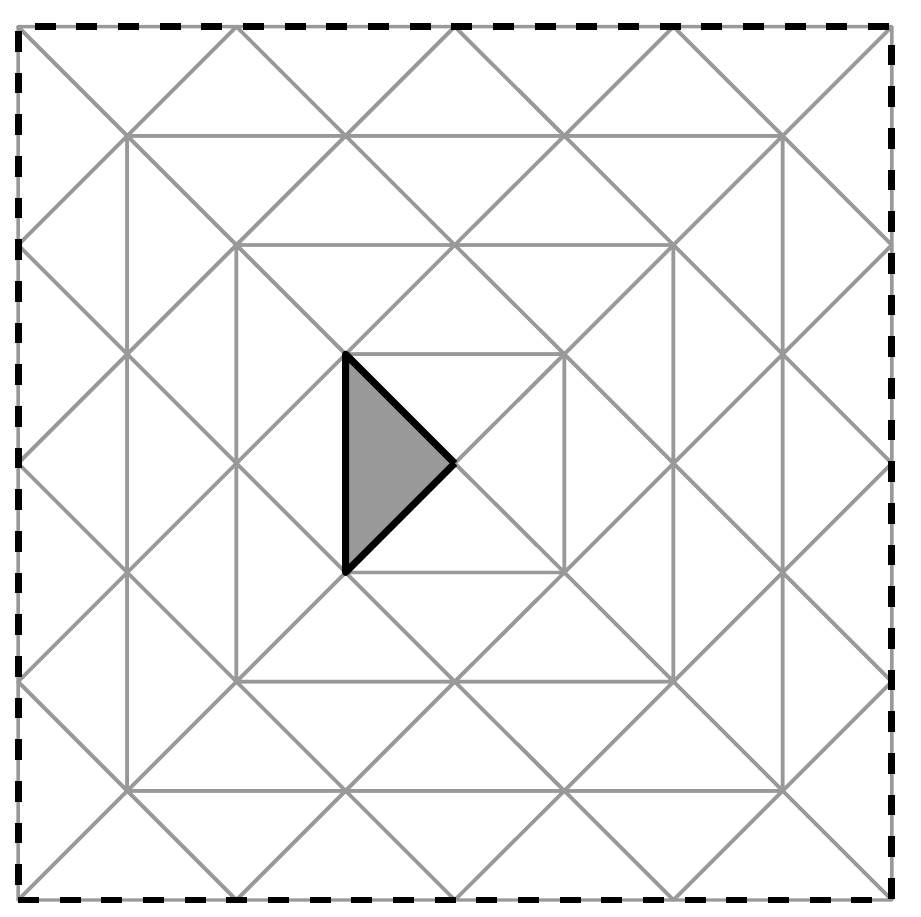}
\caption{Computational domain of the model problem of Section~\ref{s:helmholtz} and coarsest mesh.\label{fig:domain}}
\end{figure}

For the multiscale Petrov-Galerkin method, the result means that pollution effects do not occur. Note that the resolution condition $H\kappa \leq c_{\operatorname{res}}$ is somewhat minimal, because any meaningful approximation of the highly oscillatory solution of \eqref{e:modelhelm} requires at least $5-10$ degrees of freedom per wave length and coordinate direction. Saying this, we assume that the fine scale corrector problems are solved sufficiently accurate; see \cite{Gallistl.Peterseim:2015,Peterseim2014} for details.  
\smallskip

We shall present a numerical experiment taken from \cite{Peterseim2014} where this version of the method was already considered experimentally. 
Consider the scattering from sound-soft scatterer occupying the triangle $\Omega_D$. The Sommerfeld radiation condition of the scattered wave is approximated by the Robin boundary condition on the boundary $\Gamma_R:=\partial\Omega_R$ of the unit square so that $\Omega:=(0,1)^2\setminus\Omega_D$ is the computational domain; see Fig.~\ref{fig:domain}. Given the wave number $\k=2^7$, the incident wave $\displaystyle u_{inc}(x):=\exp(i\k\;x\cdot \left[\cos(0.5),\sin(0.5)\right]^T)$ is prescribed via an inhomogeneous Dirichlet boundary condition on $\Gamma_D:=\partial\Omega_D$ and the scattered wave satisfies \eqref{e:modela} with $f\equiv 0$ and the boundary conditions
\begin{align*}
  u &= -u_{inc}\quad\text{on }\Gamma_D,\\
  \nabla u\cdot \nu - i\k u &= 0\quad\text{on }\Gamma_R.
\end{align*}
The error analysis extends to this setting in a straight-forward way. 
\begin{figure}[tb]
\sidecaption[t]
\includegraphics[width=0.6\textwidth]{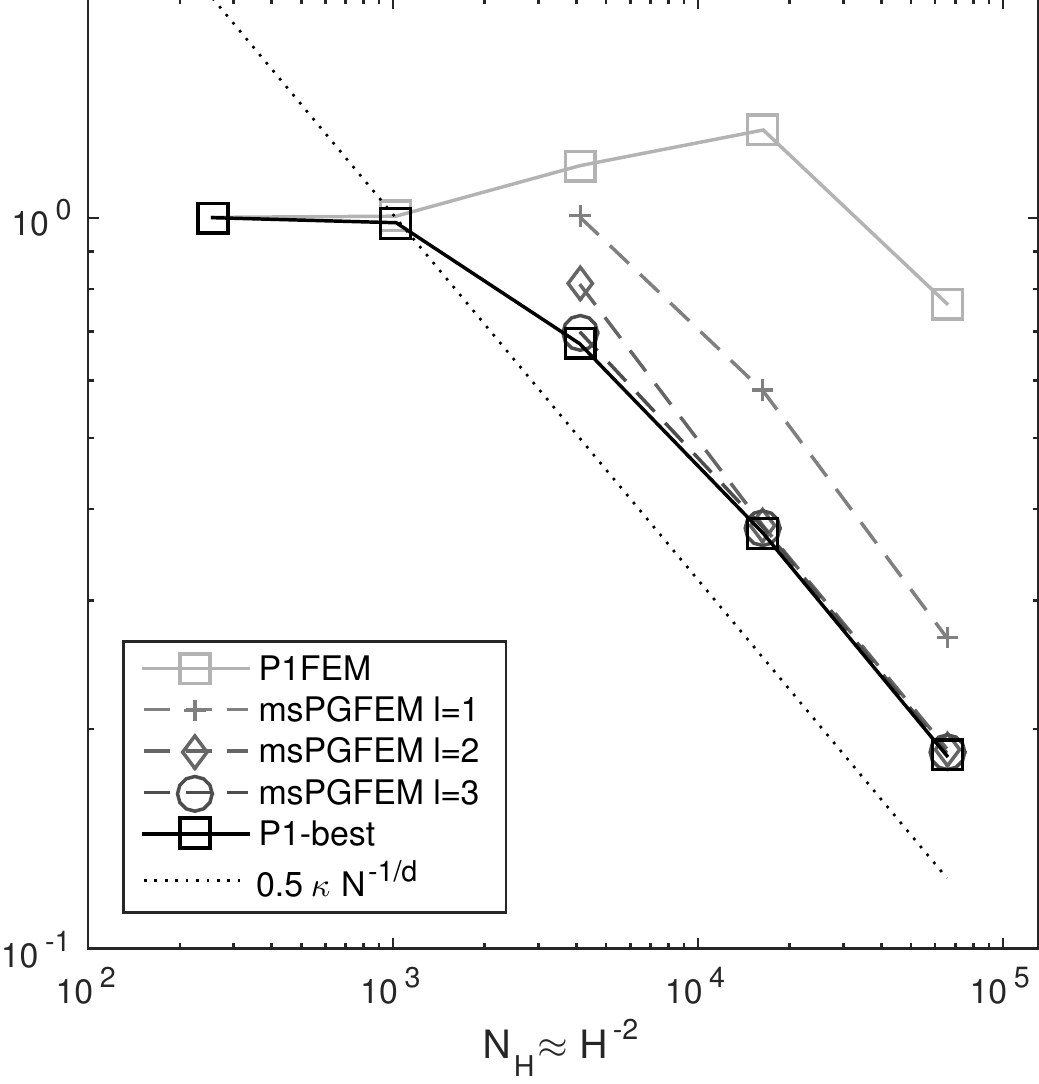}
\caption{Numerical experiment of Section~\ref{s:helmholtz}: Relative $V$-norm errors of multiscale Petrov-Galerkin method \eqref{e:discreteproblemhelm} with wave number $\k=2^7$ depending on the number of degrees of freedom $N_{H}\approx H^{-2}$, where $H=2^{-5},\ldots,2^{-7}$ is the uniform coarse mesh size. The reference mesh size $h=2^{-9}$ remains fixed. The oversampling parameter $\ell$ varies between $1$ and $3$. The $P_1$-FE solution and the best-approximation in the $P_1$-FE space on the same coarse meshes are depicted for comparison.}
 \label{fig:numhelm}
\end{figure}

We choose uniform coarse meshes of widths $H=2^{-3},\ldots,2^{-7}$ as depicted in Fig.~\ref{fig:domain}. The reference mesh $\G_h$ is derived by uniform mesh refinement of the coarse meshes and has mesh width $h=2^{-9}$. The corresponding $P_1$ conforming finite element approximation on the reference mesh $\G_h$ is denoted by $V_h$. As in the previous section, we compare the coarse scale approximations $u_{H,\ell,h}\in V_H$ with some reference solution $u_h\in V_h$.

Fig.~\ref{fig:numhelm} depicts the results for the multiscale Petrov-Galerkin method and shows that the pollution effect that is present in the $P_1$ FEM is eliminated when $\ell$ is moderately increased. For the present wave number $\ell=2$ is sufficient.  

Further numerical experiments are reported in \cite{Peterseim2014} and \cite{Gallistl.Peterseim:2015}. It is worth noting that the latter work also exploits the homogeneous structure of the PDE coefficients in the sense that only very few of the fine-scale corrector problems are actually solved due to translation invariance and symmetry. This makes the approach competitive.  
\smallskip

A very natural and straight forward generalisation of the method would be the case of heterogeneous media. The previous section plus the analysis of this section strongly indicate the potential of the method to treat high oscillations or jumps in the PDE coefficients and the pollution effect in one stroke.

\section{Final remarks}\label{s:final}
We have presented an abstract framework for the stabilization of numerical methods for multiscale partial differential equations with some focus on highly oscillatory problems. The methodology is based on the variational multiscale method and the more recent development of localised orthogonal decompositions. We have provided an abstract numerical analysis of the method which is applied to two representative model problems, a homogenization problem and a scattering problem. We have shown that the methodology can indeed eliminate critical scale-dependent pre-asymptotic effects in these cases. While the framework has already been applied successfully to other problem classes such as linear and non-linear eigenvalue problems, we expect that the framework will also be useful for convection-dominated flow, the problem that the variational method was initially designed for. 

The multiscale method presented in this paper is shown to be stable and accurate under moderate assumptions on the discretization parameters relative to characteristic parameters and length scales of the problem. These valuable properties require the pre-computation of the test basis on subgrids. These pre-computations are both local and independent, but the worst-case (serial) complexity of the method can exceed the cost of a direct numerical simulation on a global sufficiently fine mesh. If the inherent parallelism of the local cell problems cannot be exploited during the computation, we still expect a significant gain with respect to computational complexity if the pre-computation can be reused several times in the context of parameter studies, coupled problems, optimal control problems or inverse problems.
In many cases, there is also a lot of redundancy in the local problems which allows one to reduce the number of local problems drastically as it is shown in \cite{Gallistl.Peterseim:2015} in the context of acoustic scattering. We expect that this technique can be generalised to far more general situations using modern techniques of model order reduction \cite{MR2430350,MR3247814}.

We may close the discussion with some rather philosophical remark regarding the stabilization of FEMs and their inter-element continuity properties. Presently, it is believed, e.g. in the context of time-harmonic wave propagation, that stability can be increased by relaxing inter-element continuity within a discontinuous Galerkin (DG) framework. The large number of variants includes the ultra
weak variational formulation \cite{MR1618464}, Trefftz methods \cite{perugia}, DPG \cite{Zitelli20112406,dpg}, or the continuous interior penalty method
\cite{Wu2014CIPFEM}. There may be some truth in this but the general impression that relaxing continuity is the only way is certainly false as one can observe from the method presented in this paper. The multiscale Petrov-Galerkin does quite the opposite. The regularity of the test functions is increased compared to standard continuous finite elements, because they are solutions of second-order elliptic problems (at least in the ideal case). In general, test functions $w_H\in W_H$ have the property that $\ddiv A\nabla w_H\in L^2(\Omega)$. In the context of the Helmholtz model problem of Section~\ref{s:helmholtz} where $A=1$ this means that $\Delta W_H\subset L^2(\Omega)$. If $\Omega$ is convex and boundary conditions are appropriate, then $W_H\subset W^{2,2}(\Omega)$ (this can be observed for one basis function in Fig.~\ref{fig:testbasisHelmholtz}). In this respect, our methodology clearly indicates that increased differentiability might as well lead to increased stability and accuracy. Similar effects have been observed for eigenvalue computations in IGA \cite{Cottrell2006} and also LOD \cite{MP12}. This shows that breaking the inter-element continuity is not at all necessary for stability. 

\begin{acknowledgement} The present work is the result of many fruitful collaborations over the past four years \cite{MP14,EGMP13,HP13,MP12,HMP13,HMP14,Henning.Mlqvist.Peterseim:2013,Peterseim.Scheichl:2014,Brown.Peterseim:2014,Gallistl.Peterseim:2015}. I would like to thank all my co-authors, in particular Axel M{\aa}lqvist, Patrick Henning, and Dietmar Gallistl. 
\end{acknowledgement}

\newcommand{\etalchar}[1]{$^{#1}$}

\end{document}